\documentclass[12pt]{article}

\makeatletter
\renewcommand{\@oddfoot}{\hfill \thepage}
\makeatother

\usepackage[T2A]{fontenc}
\usepackage[cp1251]{inputenc}
\usepackage[russian, english]{babel}
\usepackage{mathrsfs}
\usepackage{amssymb}
\usepackage{amsmath,amsthm,amsfonts}
\usepackage[dvipdf]{graphicx}
\usepackage{fancyhdr}

\begin{document}

\begin{center}
{\bf THE EXPLICIT PROBABILITY DISTRIBUTION \\ OF THE SUM OF TWO TELEGRAPH PROCESSES}
\end{center}

\begin{center}
Alexander D. KOLESNIK\\
Institute of Mathematics and Computer Science\\
Academy Street 5, Kishinev 2028, Moldova\\
E-Mail: kolesnik@math.md
\end{center}

\vskip 0.2cm

\begin{abstract}
We consider two independent Goldstein-Kac telegraph processes $X_1(t)$ and $X_2(t)$ on the real line $\Bbb R$, both developing with finite constant speed $c>0$, 
that, at the initial time instant $t=0$, simultaneously start from the origin $0\in\Bbb R$ and whose evolutions are controlled by two independent 
homogeneous Poisson processes of the same rate $\lambda>0$. Closed-form expressions for the transition density $p(x,t)$ and 
the probability distribution function $\Phi(x,t)=\text{Pr} \{ S(t)<x \}, \; x\in\Bbb R, \; t>0,$ of the sum $S(t)=X_1(t)+X_2(t)$ of these processes 
at arbitrary time instant $t>0$, are obtained. It is also proved that the shifted time derivative $g(x,t)=(\partial/\partial t+2\lambda)p(x,t)$ 
satisfies the Goldstein-Kac telegraph equation with doubled parameters $2c$ and $2\lambda$. From this fact it follows that $p(x,t)$ solves a third-order 
hyperbolic partial differential equation, but is not its fundamental solution. The general case is also discussed. 
\end{abstract}

\vskip 0.1cm


\vskip 0.2cm


\section{Introduction}

\numberwithin{equation}{section}

The problem of summation of random variables is one of the most important fields of probability theory. The classical result states that 
if $\xi_1, \; \xi_2$ are two independent random variables with given distributions, then the distribution of their sum $\xi_1+\xi_2$ is given 
by the convolution of their distributions. Clearly, this result is also valid for the sum of arbitrary finite number of independent random variables. 
The same concerns independent stochastic processes. If $X_1(t), \; X_2(t), \; t\ge 0,$ are two independent stochastic processes with given distributions 
(that is, if for any $\tau>0, \; s>0$ the random variables $X_1(\tau), \; X_2(s)$ are independent) then the distribution of their sum 
$X_1(t)+X_2(t)$ is given by the convolution of their distributions for any fixed $t>0$. 

While the convolution operation solves the problem of desribing the distribution of the sum of two independent stochastic processes, in practice it is almost 
useless when the distributions of these processes have non-trivial forms. Except the case of exponential-type distributions, the evaluation of 
such convolutions is a very difficult and often impracticable problem.  That is why those cases (very rare indeed), when the distribution of the sum of two 
processes with non-trivial distributions can be obtained in an explicit form, look like a miracle and excite great interest. 

In the present article we examine this problem when $X_1(t), \; X_2(t)$ are two independent Goldstein-Kac telegraph processes developing 
with some constant speed and driven by two independent Poisson processes of the same rate. This subject is motivated by the great theoretical 
importance and numerous fruitful applications of the telegraph processes in physics, biology, transport phenomena, financial modelling 
and other fields. 

Another motivation has a more general mathematical character. Let we have two functions $f_1(x,t), \; f_2(x,t), \; x\in\Bbb R, \; t>0,$ (classical or generalized) 
and suppose that each of them is a solution (partial or fundamental) to respective partial differential equation (the same or different ones). 
Let $f(x,t)=f_1(x,t)\ast f_2(x,t)$ be 
the convolution (in $x$) of these functions. What is the differential equation solved by function $f(x,t)$?  In the case when $f_1(x,t), \; f_2(x,t)$ are the 
probability densities of two independent stochastic processes, this question can be treated as the problem of obtaining a partial differential equation 
for the density of the sum of these processes for arbitrary fixed $t>0$. Despite great importance of such a problem for probability theory, analysis and 
mathematical physics, it is not solved so far. Even the order of such equation is unknown. 
Our analysis throws some light upon this problem for the case when $f_1(x,t)$ and $f_2(x,t)$ are the densities of two independent 
Goldstein-Kac telegraph processes. 

The classical telegraph process $X(t)$ is performed by the stochastic motion of a particle that moves
on the real line $\Bbb R$ at some constant finite speed $c$ and alternates two possible directions of motion (forward and backward)
at Poisson-distributed random instants of intensity $\lambda>0$. Such random walk was first introduced in the works of Goldstein \cite{gold}
and Kac \cite{kac} (of which the latter is a reprinting of an earlier 1956 article).  The most remarkable fact is that the transition density of
$X(t)$ is the fundamental solution to the hyperbolic telegraph equation (which is one of the classical equations of mathematical physics). Moreover, 
under increasing $c$ and $\lambda$, it transforms into the transition density of the standard Brownian motion on $\Bbb R$. Thus, the telegraph process
can be treated as a finite-velocity counterpart of the one-dimensional Brownian motion. The telegraph process $X(t)$ can also be treated in a more
general context of random evolutions (see \cite{pin}).

During last decades the Goldstein-Kac telegraph process and its numerous generalizations have become the subject of extensive
researches. Some properties of the solution space of the Goldstein-Kac telegraph equation were studied by
Bartlett \cite{bart2}. The process of one-dimensional random
motion at finite speed governed by a Poisson process with a
time-dependent parameter was considered by Kaplan \cite{kap}. The 
relationships between the Goldstein-Kac model and physical
processes, including some emerging effects of the relativity
theory, were thoroughly examined by Bartlett \cite{bart1}, Cane
\cite{cane1, cane2}. Formulas for the distributions of the
first-exit time from a given interval and of the maximum
displacement of the telegraph process were obtained by Pinsky
\cite[Section 0.5]{pin}, Foong \cite{foong1}, Masoliver and Weiss
\cite{mas1, mas2}. The behaviour of the telegraph process with
absorbing and reflecting barriers was examined by Foong and Kanno
\cite{foong2}, Ratanov \cite{rat1}. A one-dimensional stochastic
motion with an arbitrary number of velocities and governing Poisson processes was examined by Kolesnik \cite{kol1}. The
telegraph processes with random velocities were studied by
Stadje and Zacks \cite{sta}. The behaviour of the telegraph-type evolutions
in inhomogeneous environments were considered by Ratanov \cite{rat2}. 
A detailed analysis of the moment function of the telegraph process was done by Kolesnik \cite{kol2}.
Probabilistic methods of solving the Cauchy problems for the telegraph equation
were developed by Kac \cite{kac}, Kisynski \cite{kis}, Kabanov \cite{kab}, Turbin
and Samoilenko \cite{turb}. A generalization of the Goldstein-Kac
model for the case of a damped telegraph process with logistic
stationary distributions was given by Di Crescenzo and Martinucci
\cite{cres2}. A random motion with velocities alternating at
Erlang-distributed random times was studied by Di Crescenzo
\cite{cres1}. Formulas for the occupation time distributions of
the telegraph process were obtained by Bogachev and Ratanov \cite{br}. The explicit probability distribution 
of the Euclidean distance between two independent telegraph processes with arbitrary parameters was obtained by Kolesnik \cite{kol3}. 
The most important properties of the telegraph processes and their applications to financial modelling were presented in the 
recently published book by Kolesnik and Ratanov \cite{kol4}. 

To the best of author's knowledge, despite the great variety of works on the telegraph processes, the probability laws for their linear 
combinations were not studied in the literature so far. In the present article we take the first step in this important field and examine 
the sum $S(t)=X_1(t)+X_2(t)$ of two independent Goldstein-Kac telegraph processes $X_1(t)$ and $X_2(t)$, both with the same parameters $c, \; \lambda$, that,  
at the initial time instant $t=0$, simultaneously start from the origin $0\in\Bbb R$ of the real line $\Bbb R$.  
Despite a fairly complicated form of their densities involving modified Bessel functions, one managed to obtain the transition density 
and the probability distribution function of $S(t)$ in an explicit form. To avoid the convolution operation which is practically useless in this case, 
we apply the characteristic functions technique leading to Fourier and inverse Fourier transforms combined with some important properties of Bessel and 
hypergeometric functions. We also prove that the density of $S(t)$ satisfies a third-order hyperbolic partial differential equation with an operator 
representing a product of the telegraph operator and shifted time differential operator. Some remarks on the more general case of arbitrary parameters 
and start points are also given.

\section{Some Basic Properties of the Telegraph Process}

\numberwithin{equation}{section}

The telegraph stochastic process is performed by a particle that starts at the initial time instant $t=0$ from the origin $0\in\Bbb R$ of the real
line $\Bbb R$ and moves with some finite constant speed $c$. The initial direction of the motion (positive or negative) is taken on
with equal probabilities 1/2. The motion is driven by a homogeneous Poisson process of rate $\lambda>0$ as follows. As a
Poisson event occurs, the particle instantaneously takes on the opposite direction and keeps moving with the same speed $c$ until
the next Poisson event occurs, then it takes on the opposite direction again independently of its previous motion, and so on.
This random motion has first been studied by Goldstein \cite{gold} and Kac \cite{kac} and was called the {\it telegraph process}
afterwards.

Let $X(t)$ denote the particle's position on $\Bbb R$ at time $t>0$. Since the speed $c$ is finite, then,
at arbitrary time instant $t>0$, the distribution $\text{Pr}\{ X(t)\in
dx \}$ is concentrated in the close interval $[-ct, ct]$ which is the support of this distribution. The density $f(x,t),
\; x\in\Bbb R, \; t\ge 0,$ of the distribution $\text{Pr}\{
X(t)\in dx \}$ has the structure
$$f(x, t) = f_{s}(x, t) + f_{ac}(x, t),$$
where $f_{s}(x, t)$ and $f_{ac}(x, t)$ are the densities of the
singular (with respect to the Lebesgue measure on the line) and of
the absolutely continuous components of the distribution of $X(t)$, respectively.

The singular component of the distribution is, obviously,
concentrated at two terminal points $\pm ct$ of the interval
$[-ct, ct]$ and corresponds to the case when no one Poisson event
occurs till the time moment $t$ and, therefore, the particle does not
change its initial direction. Therefore, the probability of being
at arbitrary time instant $t>0$ at the terminal points $\pm ct$ is
\begin{equation}\label{prop1}
\text{Pr}\left\{ X(t) = ct \right\} = \text{Pr}\left\{ X(t)
= - ct \right\} = \frac{1}{2} \; e^{-\lambda t} .
\end{equation}
The absolutely continuous component of the distribution of $X(t)$
is concentrated in the open interval $(-ct, ct)$ and corresponds
to the case when at least one Poisson event occurs by the moment
$t$ and, therefore, the particle changes its initial direction.
The probability of this event is
\begin{equation}\label{prop2}
\text{Pr}\left\{ X(t) \in (-ct, ct) \right\} = 1 - e^{-\lambda t}.
\end{equation}

The principal result by Goldstein \cite{gold} and Kac \cite{kac}
states that the density $f = f(x,t), \; x\in [-ct, ct], \; t>0,$
of the distribution of $X(t)$ satisfies the hyperbolic partial differential equation
\begin{equation}\label{prop3}
\frac{\partial^2 f}{\partial t^2} + 2\lambda
\frac{\partial f}{\partial t} - c^2 \frac{\partial^2 f}{\partial x^2} = 0,
\end{equation}
(which is referred to as the {\it telegraph} or {\it damped wave}
equation) and can be found by solving (\ref{prop3}) with the initial
conditions
\begin{equation}\label{iprop3}
f(x,t)\vert_{t=0} = \delta(x), \qquad
\left.\frac{\partial f(x,t)}{\partial t}\right\vert_{t=0} = 0,
\end{equation}
where $\delta(x)$ is the Dirac delta-function. This means that the
transition density $f(x,t)$ of the process $X(t)$ is the
fundamental solution (i.e. the Green's function) of the telegraph
equation (\ref{prop3}).

The explicit form of the density $f(x,t)$ is given by the formula
(see, for instance, \cite[Section 0.4]{pin} or \cite[Section 2.5]{kol4}:
\begin{equation}\label{prop4}
\aligned
f(x,t) & = \frac{e^{-\lambda t}}{2} \left[ \delta(ct-x) + \delta(ct+x) \right]\\
& \quad + \frac{\lambda e^{-\lambda t}}{2c} \left[ I_0\left( \frac{\lambda}{c} \sqrt{c^2t^2-x^2} \right) +
\frac{ct}{\sqrt{c^2t^2-x^2}} I_1\left( \frac{\lambda}{c} \sqrt{c^2t^2-x^2} \right) \right] \Theta(ct-\vert x\vert),
\endaligned
\end{equation}
where $I_0(z)$ and $I_1(z)$ are the modified Bessel functions of zero and first orders, respectively (that is, the
Bessel functions with imaginary argument) with series representations 
\begin{equation}\label{pprop5}
I_0(z) = \sum_{k=0}^{\infty} \frac{1}{(k!)^2} \left( \frac{z}{2} \right)^{2k} \qquad
I_1(z) = \sum_{k=0}^{\infty} \frac{1}{k! \; (k+1)!} \left( \frac{z}{2} \right)^{2k+1} ,
\end{equation}
and $\Theta(x)$ is the Heaviside step function
\begin{equation}\label{pprop4}
\Theta(x) = \left\{ \aligned 1, \qquad  & \text{if} \; x>0,\\
                               0, \qquad & \text{if} \; x\le 0.
\endaligned \right.
\end{equation}

The first term in (\ref{prop4})
\begin{equation}\label{prop5}
f^s(x,t) = \frac{e^{-\lambda t}}{2} \left[ \delta(ct-x) +
\delta(ct+x) \right]
\end{equation}
represents the density (in the sense of generalized functions) of the singular part of the distribution
of $X(t)$ concentrated at two terminal points $\pm ct$ of the support $[-ct, ct]$, while the second term 
\begin{equation}\label{prop6}
f^{ac}(x,t) = \frac{\lambda e^{-\lambda t}}{2c}
\left[ I_0\left( \frac{\lambda}{c} \sqrt{c^2t^2-x^2} \right) +
\frac{ct}{\sqrt{c^2t^2-x^2}} I_1\left( \frac{\lambda}{c}
\sqrt{c^2t^2-x^2} \right) \right] \Theta(ct-\vert x\vert),
\end{equation}
is the density of the absolutely continuous part of the distribution of $X(t)$ concentrated in the open interval $(-ct, ct)$.

The probability distribution function of the Goldstein-Kac telegraph process $X(t)$ has the form (see \cite[Proposition 2]{kol3}):
\begin{equation}\label{sub2}
\aligned
& \text{Pr} \left\{ X(t) < x \right\} \\
& = \left\{ \aligned
0, \hskip 5cm & x\in (-\infty, -ct],\\
\frac{1}{2} + \frac{\lambda x  e^{-\lambda t}}{2c} \sum_{k=0}^{\infty} \frac{1}{(k!)^2} \left(\frac{\lambda t}{2}\right)^{2k} \left( 1+ \frac{\lambda t}{2k+2} \right) F\left( -k, \frac{1}{2}; \frac{3}{2}; \frac{x^2}{c^2t^2} \right) , \quad & x\in (-ct, ct],\\
1, \hskip 5cm & x\in (ct, +\infty) ,
\endaligned \right. \endaligned
\end{equation}
where 
\begin{equation}\label{sub3}
F(\alpha,\beta; \gamma; z) = \sum_{k=0}^{\infty} \frac{(\alpha)_k \; (\beta)_k}{(\gamma)_k} \; \frac{z^k}{k!}
\end{equation}
is the Gauss hypergeometric function. 

The characteristic function of the telegraph process starting from the origin $x=0$ with density (\ref{prop4}) 
is given by the formula (see \cite[Section 2.4]{kol4}):
\begin{equation}\label{prop10}
\aligned
H(\xi, t) = e^{-\lambda t} & \left\{ \left[
\cosh\left( t\sqrt{\lambda^2-c^2\xi^2} \right) +
\frac{\lambda}{\sqrt{\lambda^2-c^2\xi^2}} \; \sinh\left(
t\sqrt{\lambda^2-c^2\xi^2} \right) \right] \bold 1_{\left\{\vert\xi\vert\le\frac{\lambda}{c}\right\}} \right. \\
& + \left. \left[ \cos\left( t\sqrt{c^2\xi^2 - \lambda^2}
\right) + \frac{\lambda}{\sqrt{c^2\xi^2-\lambda^2}} \;
\sin\left( t\sqrt{c^2\xi^2-\lambda^2} \right)  \right] \bold 1_{\left\{\vert\xi\vert > \frac{\lambda}{c}\right\}} \right\} ,
\endaligned
\end{equation}
where $\bold 1_{ \{ z \}}$ is the indicator function, $\xi\in\Bbb R, \; t\ge 0$.

\section{Density of the Sum of Telegraph Processes} 

\numberwithin{equation}{section}

Consider two independent telegraph processes $X_1(t)$ and $X_2(t)$ on the real line $\Bbb R$. We assume that $X_1(t)$ and $X_2(t)$ start simultaneously from the origin $x=0$ at the initial time instant $t=0$ and are developing with the same constant speed $c$. The motions are controlled by two independent  Poisson processes of the same rate $\lambda>0$, as described above.

Consider the sum
$$S(t) = X_1(t) + X_2(t) , \qquad t\ge 0,$$
of these telegraph processes. The support of the distribution $\text{Pr} \{ S(t)< x \} , \; x\in\Bbb R, \; t>0,$ of the process $S(t)$ is the close interval $[-2ct, 2ct]$. This distribution consists of two components. The singular component is concentrated at three points $0, \pm 2ct$ and corresponds to the case when no one Poisson event occurs up to time $t$. If both the processes $X_1(t)$ and $X_2(t)$ initially take the same direction (the probability of this event is  $1/2$) and no one Poisson event occurs up to time $t$ then, at moment $t$ the process $S(t)$ is located at one of the terminal points $\pm 2ct$. Thus,
\begin{equation}\label{ssumm1}
\text{Pr} \{ S(t) = 2ct \} =  \text{Pr} \{ S(t) = -2ct \} = \frac{1}{4} \; e^{-2\lambda t}, \qquad t>0 .
\end{equation}
If the processes $X_1(t)$ and $X_2(t)$ initially take different directions (the probability of this event is  $1/2$) and no one Poisson event occurs up to time $t$ then, at moment $t$ the process $S(t)$ is located at the origin and therefore
\begin{equation}\label{ssummm2}
\text{Pr} \{ S(t) = 0 \} = \frac{1}{2} \; e^{-2\lambda t}, \qquad t>0 .
\end{equation}
The remaining part $M_t = (-2ct,0) \cup (0,2ct)$ of the interval $[-2ct, 2ct]$ is the support of the absolutely continuous component of the distribution 
$\text{Pr} \{ S(t)< x \} , \; x\in\Bbb R, \; t>0,$ corresponding to the case when at least one Poisson event occurs up to time instant $t$ and therefore
\begin{equation}\label{sum1}
\text{Pr} \{ S(t) \in M_t \} = 1 - e^{-2\lambda t}, \qquad t>0 .
\end{equation}

Let $p(x,t), \; x\in\Bbb R, \; t>0,$ be the density of the process $S(t)$ treated as a generalized function. Since $X_1(t)$ and $X_2(t)$ are independent, then, for any fixed $t>0$, the density of $S(t)$ is formally given by the convolution
$$p(x,t) = f(x,t) \ast f(x,t) = \int f(z,t) \; f(x-z,t) \; dz , \qquad x\in\Bbb R, \quad t>0 ,$$
where $f(x,t)$ is the density of the telegraph processes $X_1(t)$ and $X_2(t)$ given by (\ref{prop4}). However, it seems impossible to explicitly compute this convolution due to highly complicated form of density $f(x,t)$ containing modified Bessel functions. Instead, we apply another way of finding the density $p(x,t)$ based on the  characteristic functions technique and using some important properties of special functions.

The main result of this section is given by the following theorem.

\bigskip

{\bf Theorem 1.} {\it The transition probability density $p(x,t)$ of process $S(t)$ has the form:}
\begin{equation}\label{sum2}
\aligned
p(x,t) & = \frac{e^{-2\lambda t}}{2} \; \delta(x) + \frac{e^{-2\lambda t}}{4} \left[ \delta(2ct+x) + \delta(2ct-x) \right] \\
& \quad + \frac{e^{-2\lambda t}}{2c} \left[ \lambda I_0 \left( \frac{\lambda}{c} \sqrt{4c^2t^2-x^2} \right) + \frac{1}{4} \; \frac{\partial}{\partial t} I_0 \left( \frac{\lambda}{c} \sqrt{4c^2t^2-x^2} \right) \right. \\
& \qquad\qquad\qquad + \left. \frac{\lambda^2}{2c} \int_{|x|}^{2ct} I_0 \left( \frac{\lambda}{c} \sqrt{\tau^2-x^2} \right) d\tau \right] 
\Theta(2ct-|x|) ,
\endaligned
\end{equation}
$$x\in\Bbb R, \qquad t\ge 0 .$$

\vskip 0.2cm

{\it Remark 2.} In (\ref{sum2}) the term
$$p_s(x,t) = \frac{e^{-2\lambda t}}{2} \; \delta(x) + \frac{e^{-2\lambda t}}{4} \left[ \delta(2ct+x) + \delta(2ct-x) \right]$$
represents the singular part of the density concentrated at three points 0 and $\pm 2ct$. The second term of (\ref{sum2})
\begin{equation}\label{summ3}
\aligned
p_{ac}(x,t) & = \frac{e^{-2\lambda t}}{2c} \left[ \lambda I_0 \left( \frac{\lambda}{c} \sqrt{4c^2t^2-x^2} \right) + \frac{1}{4} \; \frac{\partial}{\partial t} I_0 \left( \frac{\lambda}{c} \sqrt{4c^2t^2-x^2} \right) \right. \\
& \qquad\qquad\qquad + \left. \frac{\lambda^2}{2c} \int_{|x|}^{2ct} I_0 \left( \frac{\lambda}{c} \sqrt{\tau^2-x^2} \right) d\tau \right]\Theta(2ct-|x|) 
\endaligned
\end{equation}
represents the absolutely continuous part of the density concentrated in $M_t$. 

\vskip 0.2cm

{\it Proof.} Since the processes $X_1(t)$ and $X_2(t)$ are independent, then the characteristic function of their sum $S(t)$ is
\begin{equation}\label{sum3}
\aligned
\Psi(\xi, t) & = H^2(\xi, t) \\
& = e^{-2\lambda t} \biggl\{ \biggl[ \cosh\left( t\sqrt{\lambda^2-c^2\xi^2} \right) + \frac{\lambda}{\sqrt{\lambda^2-c^2\xi^2}} \; 
\sinh\left( t\sqrt{\lambda^2-c^2\xi^2} \right) \biggr]^2 \bold 1_{\left\{\vert\xi\vert\le\frac{\lambda}{c}\right\}} \\
& \qquad + \biggl[ \cos\left( t\sqrt{c^2\xi^2 - \lambda^2} \right) + \frac{\lambda}{\sqrt{c^2\xi^2-\lambda^2}} \;
\sin\left( t\sqrt{c^2\xi^2-\lambda^2} \right)  \biggr]^2 \bold 1_{\left\{\vert\xi\vert > \frac{\lambda}{c}\right\}} \biggr\} ,
\endaligned
\end{equation}
where $H(\xi,t)$ is the characteristic function of the telegraph process given by (\ref{prop10}). Equality (\ref{sum3}) can be represented as follows:
$$\aligned
\Psi(\xi, t) & = e^{-2\lambda t} \biggl\{ \biggl[ \cosh^2\left( t\sqrt{\lambda^2-c^2\xi^2} \right) \bold 1_{\left\{\vert\xi\vert\le\frac{\lambda}{c}\right\}} + \cos^2\left( t\sqrt{c^2\xi^2 - \lambda^2} \right) \bold 1_{\left\{\vert\xi\vert > \frac{\lambda}{c}\right\}} \biggr]  \\
& \qquad \; +  \lambda \biggl[ \frac{\sinh\left( 2t\sqrt{\lambda^2-c^2\xi^2} \right)}{\sqrt{\lambda^2-c^2\xi^2}} \bold 1_{\left\{\vert\xi\vert\le\frac{\lambda}{c}\right\}} + \frac{\sin\left( 2t\sqrt{c^2\xi^2 - \lambda^2} \right)}{\sqrt{c^2\xi^2 - \lambda^2}} \bold 1_{\left\{\vert\xi\vert > \frac{\lambda}{c}\right\}} \biggr] \\
& \qquad +  \lambda^2 \biggl[ \frac{\sinh^2\left( t\sqrt{\lambda^2-c^2\xi^2} \right)}{\lambda^2-c^2\xi^2} \bold 1_{\left\{\vert\xi\vert\le\frac{\lambda}{c}\right\}} + \frac{\sin^2\left( t\sqrt{c^2\xi^2 - \lambda^2} \right)}{c^2\xi^2 - \lambda^2} \bold 1_{\left\{\vert\xi\vert > \frac{\lambda}{c}\right\}} \biggr] \biggr\}.
\endaligned$$
Therefore, the inverse Fourier transformation of this expression yields
\begin{equation}\label{sum4}
\aligned
p(x,t) & = e^{-2\lambda t} \biggl\{ \mathcal F_{\xi}^{-1} \biggl[ \cosh^2\left( t\sqrt{\lambda^2-c^2\xi^2} \right) \bold 1_{\left\{\vert\xi\vert\le\frac{\lambda}{c}\right\}} + \cos^2\left( t\sqrt{c^2\xi^2 - \lambda^2} \right) \bold 1_{\left\{\vert\xi\vert > \frac{\lambda}{c}\right\}} \biggr](x) \\
& \quad \; +  \lambda \mathcal F_{\xi}^{-1} \biggl[ \frac{\sinh\left( 2t\sqrt{\lambda^2-c^2\xi^2} \right)}{\sqrt{\lambda^2-c^2\xi^2}} \bold 1_{\left\{\vert\xi\vert\le\frac{\lambda}{c}\right\}} + \frac{\sin\left( 2t\sqrt{c^2\xi^2 - \lambda^2} \right)}{\sqrt{c^2\xi^2 - \lambda^2}} \bold 1_{\left\{\vert\xi\vert > \frac{\lambda}{c}\right\}} \biggr](x) \\
& \quad +  \lambda^2 \mathcal F_{\xi}^{-1} \biggl[ \frac{\sinh^2\left( t\sqrt{\lambda^2-c^2\xi^2} \right)}{\lambda^2-c^2\xi^2} \bold 1_{\left\{\vert\xi\vert\le\frac{\lambda}{c}\right\}} + \frac{\sin^2\left( t\sqrt{c^2\xi^2 - \lambda^2} \right)}{c^2\xi^2 - \lambda^2} \bold 1_{\left\{\vert\xi\vert > \frac{\lambda}{c}\right\}} \biggr](x) \biggr\} .
\endaligned
\end{equation}
Our aim now is to explicitly compute inverse Fourier transforms on the right-hand side of (\ref{sum4}). For the first term in curl brackets of (\ref{sum4}) we have:

\begin{equation}\label{sum5}
\aligned
& \mathcal F_{\xi}^{-1} \left[ \cosh^2\left( t\sqrt{\lambda^2-c^2\xi^2} \right) \bold 1_{\left\{\vert\xi\vert\le\frac{\lambda}{c}\right\}} + \cos^2\left( t\sqrt{c^2\xi^2 - \lambda^2} \right) \bold 1_{\left\{\vert\xi\vert > \frac{\lambda}{c}\right\}} \right](x) \\
& = \frac{1}{2} \mathcal F_{\xi}^{-1} \left[ \left\{ \cosh\left( 2t\sqrt{\lambda^2-c^2\xi^2} \right) + 1 \right\} \bold 1_{\left\{\vert\xi\vert\le\frac{\lambda}{c}\right\}} + \left\{ \cos\left( 2t\sqrt{\lambda^2-c^2\xi^2} \right) + 1 \right\} \bold 1_{\left\{\vert\xi\vert > \frac{\lambda}{c}\right\}} \right](x) \\
& = \frac{1}{2} \delta(x) + \frac{1}{2} \mathcal F_{\xi}^{-1} \left[ \cosh\left( 2t\sqrt{\lambda^2-c^2\xi^2} \right) \bold 1_{\left\{\vert\xi\vert\le\frac{\lambda}{c}\right\}} +  \cos\left( 2t\sqrt{\lambda^2-c^2\xi^2} \right) \bold 1_{\left\{\vert\xi\vert > \frac{\lambda}{c}\right\}} \right](x) \\
& \qquad ( \text{see formula (\ref{lem6}) below}) \\
& = \frac{1}{2} \delta(x) + \frac{1}{4} \bigl[ \delta(2ct-x) + \delta(2ct+x) \bigr] + \frac{1}{8c} \; \frac{\partial}{\partial t} I_0 \left( \frac{\lambda}{c} \sqrt{4c^2t^2-x^2} \right) \Theta(2ct-|x|) .
\endaligned
\end{equation}
According to formula (\ref{lem5}) (see below), for the second term in curl brackets of (\ref{sum4}) we have:
\begin{equation}\label{sum6}
\aligned
\mathcal F_{\xi}^{-1} \biggl[ \frac{\sinh\left( 2t\sqrt{\lambda^2-c^2\xi^2} \right)}{\sqrt{\lambda^2-c^2\xi^2}} \bold 1_{\left\{\vert\xi\vert\le\frac{\lambda}{c}\right\}} + \frac{\sin\left( 2t\sqrt{c^2\xi^2 - \lambda^2} \right)}{\sqrt{c^2\xi^2 - \lambda^2}} \bold 1_{\left\{\vert\xi\vert > \frac{\lambda}{c}\right\}} \biggr](x) \\
= \frac{1}{2c} \; I_0 \left( \frac{\lambda}{c} \sqrt{4c^2t^2-x^2} \right) \Theta(2ct-|x|) .
\endaligned
\end{equation}
Finally, according to formula (\ref{lem8}) (see below), we have for the third term of (\ref{sum4}):
\begin{equation}\label{sum7}
\aligned
\mathcal F_{\xi}^{-1} & \biggl[ \frac{\sinh^2\left( t\sqrt{\lambda^2-c^2\xi^2} \right)}{\lambda^2-c^2\xi^2} \bold 1_{\left\{\vert\xi\vert\le\frac{\lambda}{c}\right\}} + \frac{\sin^2\left( t\sqrt{c^2\xi^2 - \lambda^2} \right)}{c^2\xi^2 - \lambda^2} \bold 1_{\left\{\vert\xi\vert > \frac{\lambda}{c}\right\}} \biggr](x) \\
& \qquad = \frac{1}{4c^2} \left\{ \int_{|x|}^{2ct} I_0 \left( \frac{\lambda}{c} \sqrt{\tau^2-x^2} \right) d\tau \right\} \Theta(2ct-|x|) .
\endaligned
\end{equation}
Substituting now (\ref{sum5}), (\ref{sum6}) and (\ref{sum7}) into (\ref{sum4}) we obtain (\ref{sum2}).

It remains to check that non-negative function (\ref{sum2}), being integrated in the support $[-2ct,2ct]$ of the process $S(t)$, yields 1.
Since, as is easy to see, for arbitrary $t>0$
$$\int_{-2ct}^{2ct} p_s(x,t) \; dx = e^{-2\lambda t} ,$$
then, according to (\ref{sum1}), we should verify that the absolutely continuous part $p_{ac}(x,t)$ of density (\ref{sum2}) satisfies the equality
\begin{equation}\label{sum8}
\int_{-2ct}^{2ct} p_{ac}(x,t) \; dx = 1 - e^{-2\lambda t}, \qquad t>0 .
\end{equation}
We have
\begin{equation}\label{sum9}
\aligned
\int_{-2ct}^{2ct} p_{ac}(x,t) \; dx & = \frac{e^{-2\lambda t}}{2c} \left[ \lambda \int_{-2ct}^{2ct} I_0 \left( \frac{\lambda}{c} \sqrt{4c^2t^2-x^2} \right) dx + \frac{1}{4} \int_{-2ct}^{2ct} \frac{\partial}{\partial t} I_0 \left( \frac{\lambda}{c} \sqrt{4c^2t^2-x^2} \right) dx \right. \\
& \qquad\qquad\qquad + \left. \frac{\lambda^2}{2c} \int_{-2ct}^{2ct} \left\{ \int_{|x|}^{2ct} I_0 \left( \frac{\lambda}{c} \sqrt{\tau^2-x^2} \right) d\tau \right\} dx \right] .
\endaligned
\end{equation}
According to (\ref{lemm1}), the first integral in (\ref{sum9}) is:
\begin{equation}\label{sum10}
\int_{-2ct}^{2ct} I_0 \left( \frac{\lambda}{c} \sqrt{4c^2t^2-x^2} \right) dx = \frac{2c}{\lambda} \sinh{(2\lambda t)} .
\end{equation}
Using (\ref{sum10}), we have for the second integral in (\ref{sum9}):
\begin{equation}\label{sum11}
\aligned
\int_{-2ct}^{2ct} \frac{\partial}{\partial t} I_0 \left( \frac{\lambda}{c} \sqrt{4c^2t^2-x^2} \right) dx & = \frac{\partial}{\partial t} \int_{-2ct}^{2ct} I_0 \left( \frac{\lambda}{c} \sqrt{4c^2t^2-x^2} \right) dx - 4c \\
& = \frac{\partial}{\partial t} \left( \frac{2c}{\lambda} \sinh{(2\lambda t)} \right) - 4c \\
& = 4c \left(\cosh{(2\lambda t)} - 1 \right) .
\endaligned
\end{equation}
Finally, applying formula (\ref{lem1}), we obtain for the third integral in (\ref{sum9}):
\begin{equation}\label{sum12}
\aligned
\int_{-2ct}^{2ct} & \left\{ \int_{|x|}^{2ct} I_0 \left( \frac{\lambda}{c} \sqrt{\tau^2-x^2} \right) d\tau \right\} dx \\
& = \int_{-2ct}^{2ct} \left\{ \int_0^{2ct} I_0 \left( \frac{\lambda}{c} \sqrt{\tau^2-x^2} \right) \Theta(\tau-|x|) \; d\tau \right\} dx \\ 
& = \int_0^{2ct} \left\{ \int_{-2ct}^{2ct} I_0 \left( \frac{\lambda}{c} \sqrt{\tau^2-x^2} \right) \Theta(\tau-|x|) \; dx \right\} d\tau \\
& = \int_0^{2ct} \left\{ \int_{-\tau}^{\tau} I_0 \left( \frac{\lambda}{c} \sqrt{\tau^2-x^2} \right) dx \right\} d\tau \\
& = \int_0^{2ct} \frac{2c}{\lambda} \sinh{ \left( \frac{\lambda}{c} \tau \right)} \; d\tau \\
& = \frac{2c^2}{\lambda^2} \left( \cosh{(2\lambda t)} - 1 \right) .
\endaligned
\end{equation}
Here the change of integration order is justified because the interior integral in curl brackets on the left-hand side of (\ref{sum12}) converges uniformly in $x\in (-2ct, \; 2ct)$. This fact can easily be proved by applying the mean value theorem and taking into account that $I_0(z)$ is strictly positive and monotonously  increasing continuous function. 

Substituting now (\ref{sum10}), (\ref{sum11}) and (\ref{sum12}) into (\ref{sum9}) we obtain
$$\aligned
\int\limits_{-2ct}^{2ct} p_{ac}(x,t) \; dx & = \frac{e^{-2\lambda t}}{2c} \left[ \lambda \frac{2c}{\lambda} \sinh{(2\lambda t)} + \frac{1}{4} 4c (\cosh{(2\lambda t)} - 1) + \frac{\lambda^2}{2c} \frac{2c^2}{\lambda^2} \left( \cosh{(2\lambda t)} - 1 \right) \right] \\
& = e^{-2\lambda t} \left( e^{2\lambda t} - 1 \right) \\
& = 1 - e^{-2\lambda t}
\endaligned$$
proving (\ref{sum8}). The theorem is thus completely proved. $\square$

\bigskip

The shape of the absolutely continuous part $p_{ac}(x,t)$ of the density of $S(t)$ given by (\ref{summ3}) is presented in Fig. 1.

\begin{center}
\begin{figure}[htbp]
\centerline{\includegraphics[width=10cm,height=8cm]{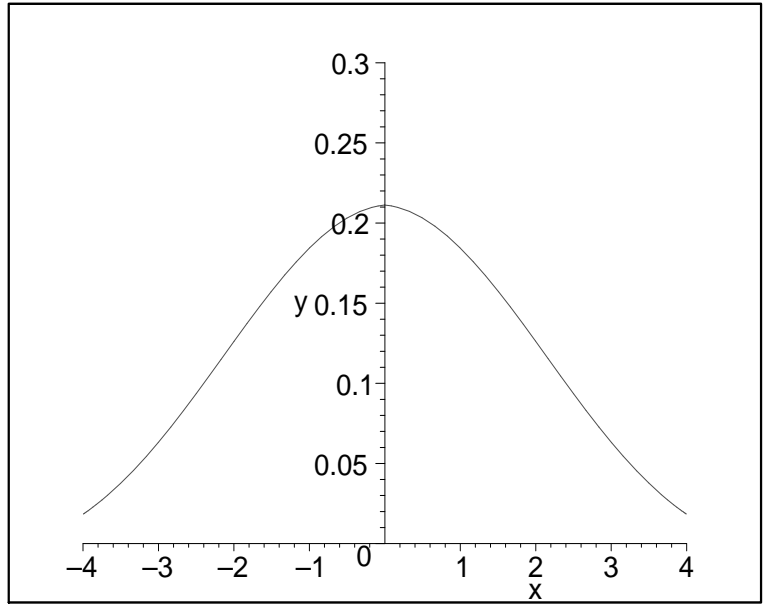}}
\caption{\it The shape of density $p_{ac}(x,t)$ at instant $t=2$ (for $c=1, \; \lambda =1$)}
\end{figure}
\end{center}

{\it Remark 3.} It is easy to check that
\begin{equation}\label{sum13}
\frac{\partial}{\partial t} I_0 \left( \frac{\lambda}{c} \sqrt{4c^2t^2-x^2} \right) = \frac{4\lambda ct}{\sqrt{4c^2t^2-x^2}} \; I_1 \left( \frac{\lambda}{c} \sqrt{4c^2t^2-x^2} \right) ,
\end{equation}
where $I_1(z)$ is the modified Bessel function of first order (see (\ref{pprop5})) and, therefore, density (\ref{sum2}) has the following alternative form: 
\begin{equation}\label{sum2a}
\aligned
p(x,t) & = \frac{e^{-2\lambda t}}{2} \; \delta(x) + \frac{e^{-2\lambda t}}{4} \left[ \delta(2ct+x) + \delta(2ct-x) \right] \\
& \quad + \frac{\lambda e^{-2\lambda t}}{2c} \left[ I_0 \left( \frac{\lambda}{c} \sqrt{4c^2t^2-x^2} \right) + \frac{ct}{\sqrt{4c^2t^2-x^2}} \; I_1 \left( \frac{\lambda}{c} \sqrt{4c^2t^2-x^2} \right) \right. \\
& \qquad\qquad\qquad + \left. \frac{\lambda}{2c} \int_{|x|}^{2ct} I_0 \left( \frac{\lambda}{c} \sqrt{\tau^2-x^2} \right) d\tau \right] 
\Theta(2ct-|x|) ,
\endaligned
\end{equation}
$$x\in\Bbb R, \qquad t\ge 0 .$$

\section{Partial Differential Equation} 

\numberwithin{equation}{section}

Consider the function
\begin{equation}\label{sum19}
g(x,t) = \left( \frac{\partial}{\partial t} + 2\lambda \right) p(x,t) .
\end{equation}
Here $\partial/\partial t$ means differentiation in $t$ of the generalized function $p(x,t)$.
The unexpected and amazing fact is that this function satisfies the Goldstein-Kac telegraph equation with doubled parameters $2c$ and $2\lambda$. 
This result is given by the following theorem.

\bigskip

{\bf Theorem 2.} {\it Function $g(x,t)$ defined by} (\ref{sum19}) {\it satisfies the telegraph equation}
\begin{equation}\label{sum20}
\left( \frac{\partial^2}{\partial t^2} + 4\lambda \frac{\partial}{\partial t} - 4c^2 \frac{\partial^2}{\partial x^2} \right) g(x,t) = 0 .
\end{equation}

\vskip 0.2cm

{\it Proof.} Introduce a new function $w(x,t)$ by the equality
$$w(x,t) = e^{2\lambda t} g(x,t) .$$
Therefore, in order to prove the theorem, we should demonstrate that function $w(x,t)$ satisfies the equation
\begin{equation}\label{sum21}
\left( \frac{\partial^2}{\partial t^2} - 4c^2 \frac{\partial^2}{\partial x^2} - 4\lambda^2 \right) w(x,t) = 0 .
\end{equation}
According to (\ref{sum19}), we have
\begin{equation}\label{sum22}
w(x,t) = e^{2\lambda t} \left( \frac{\partial}{\partial t} + 2\lambda \right) p(x,t) = \frac{\partial}{\partial t} \left( e^{2\lambda t}  p(x,t) \right) .
\end{equation}
To avoid differentiation of generalized function $w(x,t)$, instead we use the characteristic function approach. In view of (\ref{sum21}), we need to show that the characteristic function (Fourier transform) $\hat w(x,t)$ satisfies the equation
\begin{equation}\label{sum23}
\frac{\partial^2 \hat w(\xi,t)}{\partial t^2} - 4(\lambda^2 -c^2\xi^2) \; \hat w(\xi,t) = 0 .
\end{equation}
According to (\ref{sum3}), the characteristic function $\hat w(x,t)$ has the form
$$\aligned
\hat w(\xi,t) & = \frac{\partial}{\partial t} \left( e^{2\lambda t} \Psi(\xi,t) \right) \\
& = \frac{\partial}{\partial t} \biggl\{ \biggl[ \cosh\left( t\sqrt{\lambda^2-c^2\xi^2} \right) +
\frac{\lambda}{\sqrt{\lambda^2-c^2\xi^2}} \; \sinh\left(
t\sqrt{\lambda^2-c^2\xi^2} \right) \biggr]^2 \bold 1_{ \left\{\vert\xi\vert\le\frac{\lambda}{c} \right\} } \\
& \qquad + \biggl[ \cos\left( t\sqrt{c^2\xi^2 - \lambda^2} \right) + \frac{\lambda}{\sqrt{c^2\xi^2-\lambda^2}} \;
\sin\left( t\sqrt{c^2\xi^2-\lambda^2} \right)  \biggr]^2 \bold 1_{\left\{\vert\xi\vert > \frac{\lambda}{c}\right\}} \biggr\} \\
& = \frac{\partial}{\partial t} \biggl\{ \left[ \cosh^2\left( t\sqrt{\lambda^2-c^2\xi^2} \right) \bold 1_{\left\{\vert\xi\vert\le\frac{\lambda}{c}\right\}} + \cos^2\left( t\sqrt{c^2\xi^2 - \lambda^2} \right) \bold 1_{\left\{\vert\xi\vert > \frac{\lambda}{c}\right\}} \right]  \\
& \qquad \; +  \lambda \biggl[ \frac{\sinh\left( 2t\sqrt{\lambda^2-c^2\xi^2} \right)}{\sqrt{\lambda^2-c^2\xi^2}} \bold 1_{\left\{\vert\xi\vert\le\frac{\lambda}{c}\right\}} + \frac{\sin\left( 2t\sqrt{c^2\xi^2 - \lambda^2} \right)}{\sqrt{c^2\xi^2 - \lambda^2}} \bold 1_{\left\{\vert\xi\vert > \frac{\lambda}{c}\right\}} \biggr] \\
& \qquad +  \lambda^2 \biggl[ \frac{\sinh^2\left( t\sqrt{\lambda^2-c^2\xi^2} \right)}{\lambda^2-c^2\xi^2} \bold 1_{\left\{\vert\xi\vert\le\frac{\lambda}{c}\right\}} + \frac{\sin^2\left( t\sqrt{c^2\xi^2 - \lambda^2} \right)}{c^2\xi^2 - \lambda^2} \bold 1_{\left\{\vert\xi\vert > \frac{\lambda}{c}\right\}} \biggr] \biggr\}.
\endaligned$$
where $\Psi(\xi,t)$ is the characteristic function of process $S(t)$ given by (\ref{sum3}). Evaluating this expression, 
after some simple computations we arrive to the formula 
\begin{equation}\label{sum24}
\aligned
\hat w(\xi,t) & = \biggl[ \sqrt{\lambda^2-c^2\xi^2} \sinh\left( 2t\sqrt{\lambda^2-c^2\xi^2} \right) + 2\lambda \cosh\left( 2t\sqrt{\lambda^2-c^2\xi^2} \right) \\
& \qquad\qquad\qquad\qquad\qquad + \lambda^2 \; \frac{\sinh\left( 2t\sqrt{\lambda^2-c^2\xi^2} \right)}{\sqrt{\lambda^2-c^2\xi^2}} \biggr] \bold 1_{\left\{\vert\xi\vert\le\frac{\lambda}{c}\right\}} \\
& + \biggl[ - \sqrt{c^2\xi^2 - \lambda^2} \sin\left( 2t\sqrt{c^2\xi^2 - \lambda^2} \right) + 2\lambda \cos\left( 2t\sqrt{c^2\xi^2 -\lambda^2} \right) \\
& \qquad\qquad\qquad\qquad\qquad + \lambda^2 \; \frac{\sin\left( 2t\sqrt{c^2\xi^2-\lambda^2} \right)}{\sqrt{c^2\xi^2-\lambda^2}} \biggr]  \bold 1_{\left\{\vert\xi\vert > \frac{\lambda}{c}\right\}} .
\endaligned
\end{equation}
Thus, we should prove that function (\ref{sum24}) satisfies equation (\ref{sum23}). For the first term of (\ref{sum24}) we have
$$\aligned
\frac{\partial^2}{\partial t^2} & \biggl[ \sqrt{\lambda^2-c^2\xi^2} \sinh\left( 2t\sqrt{\lambda^2-c^2\xi^2} \right) + 2\lambda \cosh\left( 2t\sqrt{\lambda^2-c^2\xi^2} \right) \\
& \qquad\qquad\qquad\qquad\qquad + \lambda^2 \; \frac{\sinh\left( 2t\sqrt{\lambda^2-c^2\xi^2} \right)}{\sqrt{\lambda^2-c^2\xi^2}} \biggr] \bold 1_{\left\{\vert\xi\vert\le\frac{\lambda}{c}\right\}} \\
& = \biggl[ 4 (\lambda^2-c^2\xi^2)^{3/2} \sinh\left( 2t\sqrt{\lambda^2-c^2\xi^2} \right) + 8\lambda (\lambda^2-c^2\xi^2) \cosh\left( 2t\sqrt{\lambda^2-c^2\xi^2} \right) \\
& \qquad\qquad\qquad\qquad\qquad + 4\lambda^2 (\lambda^2-c^2\xi^2)^{1/2} \sinh\left( 2t\sqrt{\lambda^2-c^2\xi^2} \right) \biggr] \bold 1_{\left\{\vert\xi\vert\le\frac{\lambda}{c}\right\}} \\
\endaligned$$
and therefore, for $\vert\xi\vert\le\frac{\lambda}{c}$, we obtain
$$\aligned
\frac{\partial^2 \hat w(\xi,t)}{\partial t^2} & - 4(\lambda^2-c^2\xi^2) \; \hat w(\xi,t) \\
& = \biggl[ 4(\lambda^2-c^2\xi^2)^{3/2} \sinh\left( 2t\sqrt{\lambda^2-c^2\xi^2} \right) + 8\lambda (\lambda^2-c^2\xi^2) \cosh\left( 2t\sqrt{\lambda^2-c^2\xi^2} \right) \\
& \qquad\qquad\qquad\qquad\qquad + 4\lambda^2 (\lambda^2-c^2\xi^2)^{1/2} \sinh\left( 2t\sqrt{\lambda^2-c^2\xi^2} \right) \biggr] \bold 1_{\left\{\vert\xi\vert\le\frac{\lambda}{c}\right\}} \\
& \qquad - \biggl[ 4(\lambda^2-c^2\xi^2) \biggl\{ \sqrt{\lambda^2-c^2\xi^2} \sinh\left( 2t\sqrt{\lambda^2-c^2\xi^2} \right) + 2\lambda \cosh\left( 2t\sqrt{\lambda^2-c^2\xi^2} \right) \\
& \qquad\qquad\qquad\qquad\qquad + \lambda^2 \; \frac{\sinh\left( 2t\sqrt{\lambda^2-c^2\xi^2} \right)}{\sqrt{\lambda^2-c^2\xi^2}} \biggr\} \biggr] \bold 1_{\left\{\vert\xi\vert\le\frac{\lambda}{c}\right\}} \\
& = 0
\endaligned$$
proving (\ref{sum23}). The proof for the second term of (\ref{sum24}) for $\vert\xi\vert > \frac{\lambda}{c}$ is similar. The theorem is proved. $\square$

\bigskip

{\it Remark 4.} From (\ref{sum19}) and (\ref{sum20}) it follows that the transition probability density $p(x,t)$ of process $S(t)$ satisfies 
the third-order hyperbolic partial differential equation
\begin{equation}\label{sum25}
\left( \frac{\partial}{\partial t} + 2\lambda \right) \left( \frac{\partial^2}{\partial t^2} + 4\lambda \frac{\partial}{\partial t} - 4c^2 \frac{\partial^2}{\partial x^2} \right)  p(x,t) = 0 .
\end{equation}
Note that differential operator in (\ref{sum25}) represents the product of the standard Goldstein-Kac telegraph operator with doubled parameters 
$2c, \; 2\lambda$ and shifted time differential operator.  
This interesting fact means that, while the densities of two independent telegraph processes $X_1(t)$ and $X_2(t)$ satisfy the second-order telegraph equation (\ref{prop3}), their convolution (that is, the density $p(x,t)$ of the sum $S(t)=X_1(t)+X_2(t)$) satisfies third-order equation (\ref{sum25}). By differentiating 
in $t$ the characteristic function $\Psi(\xi,t)$ given by (\ref{sum3}) one can easily show that 
$$\Psi(\xi,t)\bigr|_{t=0} = 1 , \qquad \frac{\partial \Psi(\xi,t)}{\partial t}\biggr|_{t=0} = 0 , \qquad 
\frac{\partial^2 \Psi(\xi,t)}{\partial t^2}\biggr|_{t=0} = -2c^2\xi^2 ,$$
and, therefore, in contrast to (\ref{prop4}), the density $p(x,t)$ of process $S(t)$ is not the fundamental solution to equation (\ref{sum25}).

\section{Probability Distribution Function} 

\numberwithin{equation}{section}

In this section we concentrate our efforts on deriving a closed-form expression for the probability distribution function 
$$\Phi(x,t) = \text{Pr} \left\{ S(t) < x \right\} , \qquad x\in\Bbb R, \quad t>0,$$
of the process $S(t)$. This result is given by the following theorem.

\bigskip

{\bf Theorem 3.} {\it The probability distribution function $\Phi(x,t)$ has the form:}
\begin{equation}\label{dist1}
\Phi(x,t) = \left\{ \aligned
0, \qquad\qquad & \text{if} \;\; x\in (-\infty, \; -2ct],\\
G^-(x,t) , \qquad & \text{if} \;\; x\in (-2ct, \; 0], \\
G^+(x,t) , \qquad & \text{if} \;\; x\in (0, \; 2ct], \\
1, \qquad\qquad & \text{if} \;\; x\in (2ct, \; +\infty) ,
\endaligned \right. \qquad t>0,
\end{equation}
{\it where functions $G^{\pm}(x,t)$ are given by the formula:}
\begin{equation}\label{dist2}
\aligned
G^{\pm}(x,t) & = \frac{1}{2} \pm \frac{e^{-2\lambda t}}{4} \cos\left( \frac{\lambda x}{c} \right) + \frac{\lambda x e^{-2\lambda t}}{2c} \biggl[ \sum_{k=0}^{\infty} \frac{(\lambda t)^{2k}}{(k!)^2} \left( 1 + \frac{\lambda t}{2k+2} \right)  F\left( -k, \frac{1}{2}; \frac{3}{2}; \frac{x^2}{4c^2t^2} \right) \\
& \hskip 3cm + \sum_{k=0}^{\infty} \frac{(\lambda t)^{2k+1}}{(k!)^2 \; (2k+1)} \;  _3F_2 \left( -k, -k - \frac{1}{2}, \frac{1}{2}; \; -k + \frac{1}{2}, \frac{3}{2}; \; \frac{x^2}{4c^2t^2} \right) \biggr] .
\endaligned
\end{equation}
{\it Here $F(\alpha,\beta; \gamma; z)$ is the Gauss hypergeometric function given by} (\ref{sub3}) {\it and}
\begin{equation}\label{dist3}
_3F_2(\alpha,\beta,\gamma; \xi,\zeta; z) = \sum_{k=0}^{\infty} \frac{(\alpha)_k \; (\beta)_k \; (\gamma)_k}{(\xi)_k \; (\zeta)_k} \; \frac{z^k}{k!} 
\end{equation}
{\it is the generalized hypergeometric functions.}

\vskip 0.1cm

{\it Proof.} Formula (\ref{dist1}) in the intervals $x\in (-\infty, \; -2ct]$ and $x\in (2ct, \; +\infty)$ is obvious. Therefore, it remains to prove 
(\ref{dist1}) for $x\in (-2ct, \; 2ct]$. 

Since $x=0$ is a singularity point, then for arbitrary $x\in (-2ct, \; 2ct]$ we have  
$$\Phi(x,t) = \text{Pr} \left\{ S(t) = -2ct \right\} + \text{Pr} \left\{ S(t) = 0 \right\} \Theta(x) 
+ \text{Pr} \left\{ S(t) \in R_x \right\} , $$ 
where 
$$R_x = \left\{ \aligned (-2ct, \; x), \qquad & \text{if} \; x\in (-2ct, \; 0], \\ 
                           (-2ct, \; x)-\{ 0 \}, \qquad & \text{if} \; x\in (0, \; 2ct] \endaligned \right. 
$$
and $\Theta(x)$ is the Heaviside step function given by  (\ref{pprop4}).     

Taking into account (\ref{ssumm1}) and (\ref{ssummm2}), we get 
\begin{equation}\label{case1}
\Phi(x,t) = \frac{e^{-2\lambda t}}{4} + \frac{e^{-2\lambda t}}{2} \Theta(x) + \text{Pr} \left\{ S(t) \in R_x \right\} .
\end{equation}
Thus, our aim is to evaluate the term $\text{Pr} \left\{ S(t) \in R_x \right\}$ for $x\in (-2ct, \; 2ct]$. 

Integrating the absolutely continuous part of density (\ref{sum2a}), we have for arbitrary $x\in (-2ct, \; 2ct]$:  
\begin{equation}\label{dist4}
\aligned
\text{Pr} \{ S(t)\in R_x \} & = \frac{\lambda e^{-2\lambda t}}{2c} \biggl[ \int_{-2ct}^x I_0 \left( \frac{\lambda}{c} \sqrt{4c^2t^2-z^2} \right) dz
+ ct \int_{-2ct}^x \frac{I_1 \left( \frac{\lambda}{c} \sqrt{4c^2t^2-z^2} \right)}{\sqrt{4c^2t^2-z^2}} dz \\
& \qquad\qquad\qquad  + \frac{\lambda}{2c} \int_{-2ct}^x \biggl\{ \int_{|z|}^{2ct} I_0 \left( \frac{\lambda}{c} \sqrt{\tau^2-z^2} \biggr) 
d\tau \right\} \; dz \biggr] . 
\endaligned
\end{equation} 
To evaluate the integrals on the right-hand side of (\ref{dist4}), we need the following relations (see \cite[formulas (6.3) and (6.4) therein]{kol3}):   
\begin{equation}\label{dist5}
\int I_0(b\sqrt{a^2-z^2}) \; dz
= z \sum_{k=0}^{\infty} \frac{1}{(k!)^2} \left( \frac{ab}{2} \right)^{2k} F\left( -k, \frac{1}{2}; \frac{3}{2}; \frac{z^2}{a^2} \right) + \psi_1, 
\end{equation}
\begin{equation}\label{dist6}
\int \frac{I_1(b\sqrt{a^2-z^2})}{\sqrt{a^2-z^2}} \; dz = \frac{z}{a} \sum_{k=0}^{\infty} \frac{1}{k! \; (k+1)!} \left( \frac{ab}{2} \right)^{2k+1} F\left( -k, \frac{1}{2}; \frac{3}{2}; \frac{z^2}{a^2} \right) + \psi_2, 
\end{equation}
$$|z| \le a, \quad a>0, \quad b\ge 0,$$
where $F(\alpha,\beta; \gamma; z)$ is the Gauss hypergeometric function with series representation given by the first formula of (\ref{dist3}) 
and $\psi_1, \; \psi_2$ are arbitrary functions not depending on $z$. Applying formula (\ref{dist5}) to the first integral in (\ref{dist4}), we get 
$$\int\limits_{-2ct}^x I_0 \left( \frac{\lambda}{c} \sqrt{4c^2t^2-z^2} \right) dz = x \sum_{k=0}^{\infty} \frac{(\lambda t)^{2k}}{(k!)^2} F\left( -k, \frac{1}{2};  \frac{3}{2}; \frac{x^2}{4c^2t^2} \right) + 2ct \sum_{k=0}^{\infty} \frac{(\lambda t)^{2k}}{(k!)^2} F\left( -k, \frac{1}{2};  \frac{3}{2}; 1 \right) . $$
In view of the formula  
\begin{equation}\label{dist7}
 F\left( -k, \frac{1}{2};  \frac{3}{2}; 1 \right) = \frac{(2k)!!}{(2k+1)!!} = \frac{2^k \; k!}{(2k+1)!!} , \qquad k\ge 0, 
\end{equation}
the second term is found to be 
$$2ct \sum_{k=0}^{\infty} \frac{(\lambda t)^{2k}}{(k!)^2} F\left( -k, \frac{1}{2};  \frac{3}{2}; 1 \right) = \frac{c}{\lambda} \sinh(2\lambda t) $$
and we obtain for arbitrary $x\in (-2ct, 2ct]$: 
\begin{equation}\label{dist8}
\int\limits_{-2ct}^x I_0 \left( \frac{\lambda}{c} \sqrt{4c^2t^2-z^2} \right) dz = x \sum_{k=0}^{\infty} \frac{(\lambda t)^{2k}}{(k!)^2} F\left( -k, \frac{1}{2};  \frac{3}{2}; \frac{x^2}{4c^2t^2} \right) + \frac{c}{\lambda} \sinh(2\lambda t) .
\end{equation}

According to (\ref{dist6}), the second integral in (\ref{dist4}) is 
$$\aligned 
\int_{-2ct}^x \frac{I_1 \left( \frac{\lambda}{c} \sqrt{4c^2t^2-z^2} \right)}{\sqrt{4c^2t^2-z^2}} dz & = \frac{x}{2ct} \sum_{k=0}^{\infty} \frac{(\lambda t)^{2k+1}}{k! \; (k+1)!} F\left( -k, \frac{1}{2}; \frac{3}{2}; \frac{x^2}{4c^2t^2} \right) \\ 
& \qquad\qquad + \sum_{k=0}^{\infty} \frac{(\lambda t)^{2k+1}}{k! \; (k+1)!} F\left( -k, \frac{1}{2}; \frac{3}{2}; 1 \right) .
\endaligned$$  
Applying (\ref{dist7}) one can easily show that the second term is 
$$\sum_{k=0}^{\infty} \frac{(\lambda t)^{2k+1}}{k! \; (k+1)!} F\left( -k, \frac{1}{2}; \frac{3}{2}; 1 \right) = \frac{1}{2\lambda t} \left( \cosh(2\lambda t) - 1 \right)$$
and, therefore, we obtain for arbitrary $x\in (-2ct, 2ct]$: 
\begin{equation}\label{dist9}
\int_{-2ct}^x \frac{I_1 \left( \frac{\lambda}{c} \sqrt{4c^2t^2-z^2} \right)}{\sqrt{4c^2t^2-z^2}} dz = \frac{x}{2ct} \sum_{k=0}^{\infty} \frac{(\lambda t)^{2k+1}}{k! \; (k+1)!} F\left( -k, \frac{1}{2}; \frac{3}{2}; \frac{x^2}{4c^2t^2} \right) + \frac{1}{2\lambda t} \left( \cosh(2\lambda t) - 1 \right) . 
\end{equation}

For the third (double) integral in (\ref{dist4}) we have for arbitrary $x\in (-2ct, \; 2ct]$:  
\begin{equation}\label{dist3rdInt}
\aligned 
\int_{-2ct}^x & \biggl\{ \int_{|z|}^{2ct} I_0 \left( \frac{\lambda}{c} \sqrt{\tau^2-z^2} \biggr) d\tau \right\} dz \\
& = \int_{-2ct}^x \biggl\{ \int_0^{2ct} I_0 \left( \frac{\lambda}{c} \sqrt{\tau^2-z^2} \right) \bold 1_{ \{ \tau > |z| \} } \; d\tau  \biggr\} dz \\ 
& = \int_0^{2ct} \biggl\{ \int_{-2ct}^x I_0 \left( \frac{\lambda}{c} \sqrt{\tau^2-z^2} \right) \bold 1_{ \{ |z| < \tau \} } \; dz  \biggr\} d\tau \\ 
& = \int_0^{2ct} \biggl\{ \int_{-\tau}^{\min\{x, \tau \}} I_0 \left( \frac{\lambda}{c} \sqrt{\tau^2-z^2} \right) dz  \biggr\} 
\bold 1_{ \{ -\tau < \min\{x, \tau \} \} } \; d\tau . 
\endaligned
\end{equation}
Now we should separately consider two possible cases when $x$ is negative and positive. 

\vskip 0.2cm 

$\bullet$ {\it The case $x\in (-2ct, \; 0]$}. In this case $x$ is non-positive and, therefore, (\ref{dist3rdInt}) takes the form:  
\begin{equation}\label{dist10}
\aligned 
\int_{-2ct}^x \biggl\{ \int_{|z|}^{2ct} I_0 \left( \frac{\lambda}{c} \sqrt{\tau^2-z^2} \biggr) d\tau \right\} dz & = \int_0^{2ct} \biggl\{ \int_{-\tau}^x I_0 \left( \frac{\lambda}{c} \sqrt{\tau^2-z^2} \right) dz  \biggr\} \bold 1_{ \{ -\tau < x \} } \; d\tau \\
& = \int_{-x}^{2ct} \biggl\{ \int_{-\tau}^x I_0 \left( \frac{\lambda}{c} \sqrt{\tau^2-z^2} \right) dz  \biggr\} d\tau . 
\endaligned
\end{equation}
According to (\ref{dist5}), the interior integral in curl brackets is: 
$$\aligned
\int_{-\tau}^x I_0 \left( \frac{\lambda}{c} \sqrt{\tau^2-z^2} \right) dz & = x \sum_{k=0}^{\infty} \frac{1}{(k!)^2} \left( \frac{\lambda}{2c} \tau \right)^{2k} F\left( -k, \frac{1}{2}; \frac{3}{2}; \frac{x^2}{\tau^2} \right) \\ 
& \hskip 0.5cm + \tau \sum_{k=0}^{\infty} \frac{1}{(k!)^2} \left( \frac{\lambda}{2c} \tau \right)^{2k} F\left( -k, \frac{1}{2}; \frac{3}{2}; 1 \right) .
\endaligned$$
Applying again (\ref{dist7}) we easily evaluate the second term 
$$\tau \sum_{k=0}^{\infty} \frac{1}{(k!)^2} \left( \frac{\lambda}{2c} \tau \right)^{2k} F\left( -k, \frac{1}{2}; \frac{3}{2}; 1 \right) = \frac{c}{\lambda} \sinh\left(\frac{\lambda}{c} \tau \right) $$ 
and, therefore, we get 
\begin{equation}\label{dist11}
\int_{-\tau}^x I_0 \left( \frac{\lambda}{c} \sqrt{\tau^2-z^2} \right) dz = x \sum_{k=0}^{\infty} \frac{1}{(k!)^2} \left( \frac{\lambda}{2c} \tau \right)^{2k} F\left( -k, \frac{1}{2}; \frac{3}{2}; \frac{x^2}{\tau^2} \right) + \frac{c}{\lambda} \sinh\left(\frac{\lambda}{c} \tau \right) .
\end{equation}
Substituting this into (\ref{dist10}) we obtain: 
\begin{equation}\label{dist12}
\aligned 
& \int_{-2ct}^x \biggl\{ \int_{|z|}^{2ct} I_0 \left( \frac{\lambda}{c} \sqrt{\tau^2-z^2} \biggr) d\tau \right\} \; dz \\ 
& = x \sum_{k=0}^{\infty} \frac{1}{(k!)^2} \left( \frac{\lambda}{2c} \right)^{2k} \int_{-x}^{2ct} \tau^{2k} F\left( -k, \frac{1}{2}; \frac{3}{2}; \frac{x^2}{\tau^2} \right)  d\tau + \frac{c}{\lambda} \int_{-x}^{2ct} \sinh\left(\frac{\lambda}{c} \tau \right)  d\tau .
\endaligned
\end{equation}
Here we have used the easily checked fact that, for arbitrary $x\in (-2ct, \; 2ct]$, the following uniform (in $\tau\in (0, \; 2ct)$) estimate 
$$\left| \tau^{2k} F\left( -k, \frac{1}{2}; \frac{3}{2}; \frac{x^2}{\tau^2} \right) \right| < (4ct)^{2k} $$
holds and, therefore, the series in (\ref{dist11}) converges uniformly in $\tau\in (0, \; 2ct)$. 

The second term on the right-hand side of (\ref{dist12}) is found to be 
\begin{equation}\label{dist13}
\frac{c}{\lambda} \int_{-x}^{2ct} \sinh\left(\frac{\lambda}{c} \tau \right)  d\tau = 
\frac{c^2}{\lambda^2} \left[ \cosh(2\lambda t) - \cosh\left( \frac{\lambda}{c} x \right) \right] .
\end{equation}

Applying formula (\ref{hyp1}) of Lemma A4 (see below), we have for the integral of the first term on the right-hand side of (\ref{dist12}): 
\begin{equation}\label{dist14}
\aligned 
\int_{-x}^{2ct} \tau^{2k} F\left( -k, \frac{1}{2}; \frac{3}{2}; \frac{x^2}{\tau^2} \right) d\tau 
& = \frac{(2ct)^{2k+1}}{2k+1} \; _3F_2 \left( -k, -k - \frac{1}{2}, \frac{1}{2}; \; -k + \frac{1}{2}, \frac{3}{2}; \; \frac{x^2}{4c^2t^2} \right) \\ 
& \qquad + \frac{x^{2k+1}}{2k+1} \; _3F_2 \left( -k, -k - \frac{1}{2}, \frac{1}{2}; \; -k + \frac{1}{2}, \frac{3}{2}; \; 1 \right) . 
\endaligned
\end{equation}
Using the well-known formulas for Pochhammer symbol 
\begin{equation}\label{poch}
\aligned 
(-k)_s & = (-1)^s \; \frac{k!}{(k-s)!} , \qquad 0\le s \le k, \quad k\ge 0, \\
\frac{(a)_s}{(a+1)_s} & = \frac{a}{a+s}, \qquad s\ge 0, \quad \quad a\in\Bbb R, 
\endaligned
\end{equation}
and some simple combinatorial relations, one can check that 
$$\aligned 
_3F_2 \left( -k, -k - \frac{1}{2}, \frac{1}{2}; \; -k + \frac{1}{2}, \frac{3}{2}; \; 1 \right) 
& = (2k+1) \sum_{s=0}^k \frac{(-1)^s}{(2s+1)(2k-2s+1)} \; \binom ks \\ 
& = \left\{ \aligned \frac{2^k \; k!}{(k+1) \; (2k-1)!!} , \qquad & \text{if $k$ is even} , \\ 
                  0 , \qquad\qquad & \text{if $k$ is odd} . 
\endaligned \right. 
\endaligned$$
Therefore, (\ref{dist14}) takes the form: 
 
\begin{equation}\label{dist15}
\aligned 
\int_{-x}^{2ct} \tau^{2k} F\left( -k, \frac{1}{2}; \frac{3}{2}; \frac{x^2}{\tau^2} \right) d\tau 
& = \frac{(2ct)^{2k+1}}{2k+1} \; _3F_2 \left( -k, -k - \frac{1}{2}, \frac{1}{2}; \; -k + \frac{1}{2}, \frac{3}{2}; \; \frac{x^2}{4c^2t^2} \right) \\ 
& \qquad + \left\{ \aligned \frac{2^k \; k! \;\; x^{2k+1}}{(k+1) \; (2k+1)!!} , \qquad & \text{if $k$ is even} , \\ 
                  0 , \qquad\qquad & \text{if $k$ is odd} . 
\endaligned \right. 
\endaligned
\end{equation}

Substituting now (\ref{dist15}) and (\ref{dist13}) into (\ref{dist12}) and taking into account that   
\begin{equation}\label{keven}
x \sum_{\substack{k=0\\ k \; \text{is even}}}^{\infty} \frac{1}{(k!)^2} \left( \frac{\lambda}{2c} \right)^{2k}  \frac{2^k \; k! \;\; x^{2k+1}}{(k+1) \; (2k+1)!!} = 
\frac{c^2}{\lambda^2} \left[ \cosh\left( \frac{\lambda}{c} x\right) - \cos\left( \frac{\lambda}{c} x\right) \right] ,
\end{equation}
we obtain, for arbitrary $x\in (-2ct, \; 0]$, the following formula: 
\begin{equation}\label{dist16}
\aligned 
\int_{-2ct}^x & \biggl\{ \int_{|z|}^{2ct} I_0 \left( \frac{\lambda}{c} \sqrt{\tau^2-z^2} \biggr) d\tau \right\} dz \\ 
& = \frac{2c x}{\lambda} \sum_{k=0}^{\infty} \frac{(\lambda t)^{2k+1}}{(k!)^2 \; (2k+1)} \;  
_3F_2 \left( -k, -k - \frac{1}{2}, \frac{1}{2}; \; -k + \frac{1}{2}, \frac{3}{2}; \; \frac{x^2}{4c^2t^2} \right) \\ 
& \qquad + \frac{c^2}{\lambda^2} \left[ \cosh(2\lambda t) - \cos\left( \frac{\lambda}{c} x\right) \right] , \qquad \text{for} \; x\in (-2ct, \; 0] . 
\endaligned
\end{equation}

Substituting (\ref{dist8}), (\ref{dist9}) and (\ref{dist16}) into (\ref{dist4}), after some simple computations, we get for arbitrary $x\in (-2ct, \; 0]$:
\begin{equation}\label{Gneg}
\aligned
\text{Pr} \{ S(t)\in R_x \} & = \frac{1}{2} - \frac{e^{-2\lambda t}}{2} \cos^2\left( \frac{\lambda}{2c} x \right) \\ 
& + \frac{\lambda x e^{-2\lambda t}}{2c} \biggl[ \sum_{k=0}^{\infty} \frac{(\lambda t)^{2k}}{(k!)^2} \left( 1 + \frac{\lambda t}{2k+2} \right)  F\left( -k, \frac{1}{2}; \; \frac{3}{2}; \; \frac{x^2}{4c^2t^2} \right) \\
& \hskip 1.5cm + \sum_{k=0}^{\infty} \frac{(\lambda t)^{2k+1}}{(k!)^2 \; (2k+1)} \;  _3F_2 \left( -k, -k - \frac{1}{2}, \frac{1}{2}; \; -k + \frac{1}{2}, \frac{3}{2}; \; \frac{x^2}{4c^2t^2} \right) \biggr] .  
\endaligned
\end{equation}
Substituting (\ref{Gneg}) into (\ref{case1}), we finally obtain function $G^-(x,t)$ defined in the interval $x\in (-2ct, \; 0]$ and given by formula (\ref{dist2}). 

\vskip 0.2cm

$\bullet$ {\it The case $x\in (0, \; 2ct]$}. In this case $x$ is strictly positive and, therefore, (\ref{dist3rdInt}) yields: 
\begin{equation}\label{dist18}
\aligned 
\int_{-2ct}^x & \biggl\{ \int_{|z|}^{2ct} I_0 \left( \frac{\lambda}{c} \sqrt{\tau^2-z^2} \biggr) d\tau \right\} dz \\
& = \int_0^{2ct} \biggl\{ \int_{-\tau}^{\min\{x, \tau \}} I_0 \left( \frac{\lambda}{c} \sqrt{\tau^2-z^2} \right) dz  \biggr\} \; d\tau \\
& = \int_0^{2ct} \biggl\{ \int_{-\tau}^x I_0 \left( \frac{\lambda}{c} \sqrt{\tau^2-z^2} \right) dz  \biggr\} \bold 1_{ \{ \tau\ge x \} } \; d\tau \\ 
& \hskip 2cm + \int_0^{2ct} \biggl\{ \int_{-\tau}^{\tau} I_0 \left( \frac{\lambda}{c} \sqrt{\tau^2-z^2} \right) dz  \biggr\} \bold 1_{ \{ \tau < x\} } \; d\tau \\
& = \int_x^{2ct} \biggl\{ \int_{-\tau}^x I_0 \left( \frac{\lambda}{c} \sqrt{\tau^2-z^2} \right) dz  \biggr\} \; d\tau \\ 
& \hskip 2cm + \int_0^x \biggl\{ \int_{-\tau}^{\tau} I_0 \left( \frac{\lambda}{c} \sqrt{\tau^2-z^2} \right) dz  \biggr\} \; d\tau .
\endaligned
\end{equation}
Applying (\ref{dist11}), we get for the first integral on the right-hand side of (\ref{dist18}):

$$\aligned 
\int_x^{2ct} & \biggl\{ \int_{-\tau}^x I_0 \left( \frac{\lambda}{c} \sqrt{\tau^2-z^2} \right) dz  \biggr\} \; d\tau \\
& = x \sum_{k=0}^{\infty} \frac{1}{(k!)^2} \left( \frac{\lambda}{2c} \right)^{2k} \int_x^{2ct} \tau^{2k} F\left( -k, \frac{1}{2}; \frac{3}{2}; \frac{x^2}{\tau^2} \right)  d\tau + \frac{c}{\lambda} \int_x^{2ct} \sinh\left(\frac{\lambda}{c} \tau \right)  d\tau .
\endaligned$$
Replacing $-x \mapsto x$ in (\ref{dist13}) and (\ref{dist15}), we get these integrals:
$$\frac{c}{\lambda} \int_x^{2ct} \sinh\left(\frac{\lambda}{c} \tau \right)  d\tau = \frac{c^2}{\lambda^2} \left[ \cosh(2\lambda t) - \cosh\left( \frac{\lambda}{c} x \right) \right] $$
$$\aligned 
\int_x^{2ct} \tau^{2k} F\left( -k, \frac{1}{2}; \frac{3}{2}; \frac{x^2}{\tau^2} \right) d\tau 
& = \frac{(2ct)^{2k+1}}{2k+1} \; _3F_2 \left( -k, -k - \frac{1}{2}, \frac{1}{2}; \; -k + \frac{1}{2}, \frac{3}{2}; \; \frac{x^2}{4c^2t^2} \right) \\ 
& \qquad - \left\{ \aligned \frac{2^k \; k! \;\; x^{2k+1}}{(k+1) \; (2k+1)!!} , \qquad & \text{if $k$ is even} , \\ 
                  0 , \qquad\qquad & \text{if $k$ is odd} , 
\endaligned \right. 
\endaligned$$
and, therefore, taking into account (\ref{keven}), after some simple calculations we arrive to the following formula valid for arbitrary $x\in (0, \; 2ct]$: 
\begin{equation}\label{dist19}
\aligned 
\int_x^{2ct} & \biggl\{ \int_{-\tau}^x I_0 \left( \frac{\lambda}{c} \sqrt{\tau^2-z^2} \right) dz  \biggr\} \; d\tau \\
& = \frac{2c x}{\lambda} \sum_{k=0}^{\infty} \frac{(\lambda t)^{2k+1}}{(k!)^2 \; (2k+1)} \;  
_3F_2 \left( -k, -k - \frac{1}{2}, \frac{1}{2}; \; -k + \frac{1}{2}, \frac{3}{2}; \; \frac{x^2}{4c^2t^2} \right) \\ 
& \qquad + \frac{c^2}{\lambda^2} \left[ \cosh(2\lambda t) -2\cosh\left( \frac{\lambda}{c} x \right) + \cos\left( \frac{\lambda}{c} x\right) \right] , 
\qquad \text{for} \; x\in (0, \; 2ct] . 
\endaligned
\end{equation}

Taking into account that
$$\int_{-\tau}^{\tau} I_0 \left( \frac{\lambda}{c} \sqrt{\tau^2-z^2} \right) dz = \frac{2c}{\lambda} \sinh\left( \frac{\lambda}{c} \tau \right) ,$$
we can easily evaluate the second integral on the right-hand side of (\ref{dist18}): 
\begin{equation}\label{dist20}
\int_0^x \biggl\{ \int_{-\tau}^{\tau} I_0 \left( \frac{\lambda}{c} \sqrt{\tau^2-z^2} \right) dz \biggr\} \; d\tau 
= \frac{2c^2}{\lambda^2} \left[ \cosh\left( \frac{\lambda}{c} x \right) - 1 \right] , 
\qquad \text{for} \; x\in (0, \; 2ct] . 
\end{equation}

Substituting (\ref{dist19}) and (\ref{dist20}) into (\ref{dist18}), we obtain: 
\begin{equation}\label{dist21}
\aligned 
\int_{-2ct}^x & \biggl\{ \int_{|z|}^{2ct} I_0 \left( \frac{\lambda}{c} \sqrt{\tau^2-z^2} \biggr) d\tau \right\} dz \\
& = \frac{2c x}{\lambda} \sum_{k=0}^{\infty} \frac{(\lambda t)^{2k+1}}{(k!)^2 \; (2k+1)} \;  
_3F_2 \left( -k, -k - \frac{1}{2}, \frac{1}{2}; \; -k + \frac{1}{2}, \frac{3}{2}; \; \frac{x^2}{4c^2t^2} \right) \\ 
& \qquad + \frac{c^2}{\lambda^2} \left[ \cosh(2\lambda t) + \cos\left( \frac{\lambda}{c} x\right) - 2 \right] , 
\qquad \text{for} \; x\in (0, \; 2ct] . 
\endaligned
\end{equation}

Substituting now (\ref{dist8}), (\ref{dist9}) and (\ref{dist21}) into (\ref{dist4}), after some simple computations, we obtain the following formula 
valid for arbitrary $x\in (0, \; 2ct]$: 
\begin{equation}\label{dist22}
\aligned
\text{Pr} \{ S(t)\in R_x \} & = \frac{1}{2} - \frac{e^{-2\lambda t}}{4}\left( 3 - \cos\left( \frac{\lambda}{c} x\right) \right) \\ 
& + \frac{\lambda x e^{-2\lambda t}}{2c} \biggl[ \sum_{k=0}^{\infty} \frac{(\lambda t)^{2k}}{(k!)^2} \left( 1 + \frac{\lambda t}{2k+2} \right)  F\left( -k, \frac{1}{2}; \; \frac{3}{2}; \; \frac{x^2}{4c^2t^2} \right) \\
& \hskip 1cm + \sum_{k=0}^{\infty} \frac{(\lambda t)^{2k+1}}{(k!)^2 \; (2k+1)} \;  _3F_2 \left( -k, -k - \frac{1}{2}, \frac{1}{2}; \; -k + \frac{1}{2}, \frac{3}{2}; \; \frac{x^2}{4c^2t^2} \right) \biggr] , 
\endaligned
\end{equation}
$$\text{for} \; x\in (0, \; 2ct].$$
Substituting (\ref{dist22}) into (\ref{case1}) we obtain function $G^+(x,t)$ defined in the interval $x\in (0, \; 2ct]$ and given by formula (\ref{dist2}). 
The theorem is thus completely proved. $\square$

\bigskip 

The shape of probability distribution function $\Phi(x,t)$ at time instant $t=1.5$ given by formulas (\ref{dist1})-(\ref{dist2}) for the particular values of 
parameters $c=2, \; \lambda =0.8$ in the interval $x\in (-6, \; 6]$ is presented in Fig. 2. 

\begin{center}
\begin{figure}[htbp]
\centerline{\includegraphics[width=10cm,height=8cm]{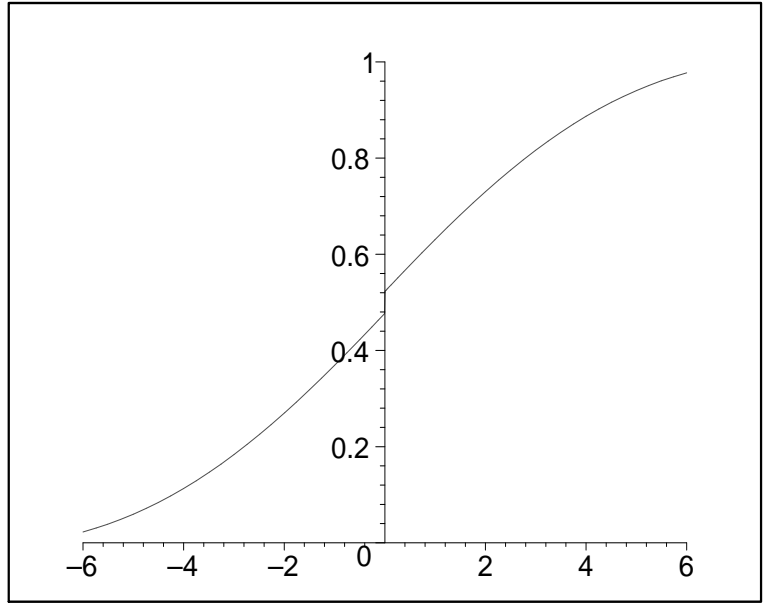}}
\caption{\it The shape of p.d.f. $\Phi(x,t)$ at time instant $t=1.5$ (for $c=2, \; \lambda =0.8$)}
\end{figure}
\end{center}

We see that distribution function $\Phi(x,t)$ is left-continuous with jumps at the origin $x=0$ 
and at the terminal points $x=\pm 2ct$ determined by the singularities concentrated at these three points. Obviously,   
$$\lim\limits_{x\to 0-0} G^-(x, t) = G^-(0, t) = \frac{1}{2} - \frac{e^{-2\lambda t}}{4} , \qquad 
\lim\limits_{x\to 0+0} G^+(x,t) = \frac{1}{2} + \frac{e^{-2\lambda t}}{4}$$
and, therefore, at the origin $x=0$ function $\Phi(x,t)$ has jump of the amplitude
$$\lim\limits_{x\to 0+0} G^+(x,t)  - \lim\limits_{x\to 0-0} G^-(x,t)  = \frac{e^{-2\lambda t}}{2} .$$
This entirely accords with (\ref{ssummm2}). One can also check that
\begin{equation}\label{dist23}
\lim\limits_{x\to -2ct+0} G^-(x,t) = \frac{e^{-2\lambda t}}{4} , \qquad \lim\limits_{x\to 2ct-0} G^+(x,t) = G^+(2ct,t) = 1 - \frac{e^{-2\lambda t}}{4}
\end{equation}
and, hence, function $\Phi(x,t)$ has jumps of the same amplitude $e^{-2\lambda t}/4$ at the terminal points $x=\pm 2ct$. This
also entirely accords with (\ref{ssumm1}). 

\bigskip

{\it Remark 5.} Notice that the first series in functions $G^{\pm}(x,t)$ containing Gauss hypergeometric function plus 1/2 is quite similar 
to the probability distribution function (\ref{sub2}) of the Goldstein-Kac telegraph process. 

\bigskip

{\it Remark 6.} Probability distribution function $\Phi(x,t)$ has very interesting and unexpected peculiarity. Both the functions $G^{\pm}(x,t)$ contain an 
oscillating term determined by the presence of the cosine function $\pm\cos\left( \frac{\lambda}{c} x\right)$. On the other hand, we see that   
$G^{\pm}(x,t)$ are strictly positive and monotonously increasing functions as they must be. The explanation of this unusual fact is that the second 
series in $G^{\pm}(x,t)$ that includes the hypergeometric function $_3F_2\left( -k, -k - \frac{1}{2}, \frac{1}{2}; \; -k + \frac{1}{2}, \frac{3}{2}; \; \frac{x^2}{4c^2t^2} \right)$ converges uniformly in $x\in [-2ct, \; 2ct]$ and contains an infinite number of terms that form a hidden oscillating term with the cosine function $\mp\cos\left( \frac{\lambda}{c} x\right)$ of the opposite sign and, therefore, these oscillating terms annihilate each other. This interesting phenomenon can be observed 
when computing the limits (\ref{dist23}). 
The appearance of such oscillating terms in functions $G^{\pm}(x,t)$ (which are, in fact, the convolution of two probability distribution functions of 
the telegraph processes $X_1(t)$ and $X_2(t)$) is a fairly unusual fact that can, apparently, be explained by some properties of the convolution operation.  

\bigskip

{\it Remark 7.} Using the relation (see \cite[item 7.4.1, formula 5]{pbm2}):
$$\aligned
_3F_2 & \left( -k, -k - \frac{1}{2}, \frac{1}{2}; \; -k + \frac{1}{2}, \frac{3}{2}; \; \frac{x^2}{4c^2t^2} \right) \\
& = \frac{1}{2k+2} F\left( -k-\frac{1}{2}, -k; \; -k + \frac{1}{2}; \; \frac{x^2}{4c^2t^2} \right) +  \frac{2k+1}{2k+2} F\left( -k, \frac{1}{2}; \; \frac{3}{2}; \; \frac{x^2}{4c^2t^2} \right) ,
\endaligned$$
we can represent probability distribution function $\Phi(x,t)$ in terms of solely Gauss hypergeometric function with the following alternative form of functions $G^{\pm}(x,t)$: 
\begin{equation}\label{dist24}
\aligned
G^{\pm}(x,t) & = \frac{1}{2} \pm \frac{e^{-2\lambda t}}{4} \cos\left( \frac{\lambda x}{c} \right) + \frac{\lambda x e^{-2\lambda t}}{2c} \biggl[ \sum_{k=0}^{\infty} \frac{(\lambda t)^{2k}}{(k!)^2} \left( 1 + \frac{\lambda t}{k+1} \right)  F\left( -k, \frac{1}{2}; \frac{3}{2}; \frac{x^2}{4c^2t^2} \right) \\
& \hskip 3cm + \sum_{k=0}^{\infty} \frac{(\lambda t)^{2k+1}}{(k!)^2 \; (2k+1) (2k+2)} \; F\left( -k, -k - \frac{1}{2}; \; -k + \frac{1}{2}; \; \frac{x^2}{4c^2t^2} \right) \biggr] .
\endaligned
\end{equation}

\section{Some Remarks on the General Case}

\numberwithin{equation}{section}

The results obtained above concern the case when both telegraph processes start from the origin $0\in\Bbb R$ and have the same parameters $c$ and $\lambda$. 
The most general case implies that the processes may have different parameters and may start from two different points of $\Bbb R$. While the method developed 
in this article works also in this situation, the analysis seems to be much more complicated and explicit formulas for the distribution of the sum of the processes 
can scarcely be obtained. In this section we give some hints concerning such general case.  

Denote by $X^{x^0}(t)$ the telegraph process starting from some arbitrary point $x^0\in\Bbb R$. It is clear that
the transition density of $X^{x^0}(t)$ emerges from (\ref{prop4}) by the formal replacement $x \mapsto x-x^0$ and it has the form 
\begin{equation}\label{gen1}
\aligned
f_{x^0}(x,t) & = \frac{e^{-\lambda t}}{2} \left[ \delta(ct-(x-x^0)) + \delta(ct+(x-x^0)) \right] \\
& + \frac{\lambda e^{-\lambda t}}{2c} \left[ I_0\left( \frac{\lambda}{c} \sqrt{c^2t^2-(x-x^0)^2} \right)
+ ct \frac{I_0\left( \frac{\lambda}{c} \sqrt{c^2t^2-(x-x^0)^2}\right)}{\sqrt{c^2t^2-(x-x^0)^2}} \right]\Theta(ct-\vert x-x^0\vert).
\endaligned
\end{equation}
The support of the distribution of $X^{x^0}(t)$ is the close interval $[x^0-ct, x^0+ct]$. The first term in (\ref{gen1})
\begin{equation}\label{gen2}
f_s^{x^0}(x,t) = \frac{e^{-\lambda t}}{2} \left[ \delta(ct-(x-x^0)) +
\delta(ct+(x-x^0)) \right]
\end{equation}
is the singular part of the density concentrated at two terminal points $x^0\pm ct$ of the interval, while the second term
\begin{equation}\label{gen3}
f_{ac}^{x^0}(x,t) = \frac{\lambda e^{-\lambda t}}{2c} \left[ I_0\left( \frac{\lambda}{c} \sqrt{c^2t^2-(x-x^0)^2} \right)
+ ct \frac{I_0\left( \frac{\lambda}{c} \sqrt{c^2t^2-(x-x^0)^2}\right)}{\sqrt{c^2t^2-(x-x^0)^2}} \right]\Theta(ct-\vert x-x^0\vert).
\end{equation}
is the density of the absolutely continuous part of the distribution of $X^{x^0}(t)$ concentrated in the open interval $(x^0-ct, x^0+ct)$.

The characteristic function of process $X^{x^0}(t)$ has the form 
\begin{equation}\label{gen4}
H^{x^0}(\xi,t) = e^{i\xi x^0} H(\xi,t) , \qquad \xi\in\Bbb R, \quad t\ge 0,
\end{equation}
where $H(\xi,t)$ is the characteristic function of the telegraph process starting from the origin $x=0$ and given by (\ref{prop10}). 
Obviously, $H^{x^0}(\xi,t)$ is a complex function if $x^0\neq 0$.  

Let $X_1^{x_1^0}(t)$ and $X_2^{x_2^0}(t)$ be two independent telegraph processes that, at the initial time instant $t=0$, simultaneously 
start from two arbitrary points $x_1^0, x_2^0\in\Bbb R$, respectively. The general case implies that $X_1^{x_1^0}(t)$ and $X_2^{x_2^0}(t)$ 
develop with arbitrary constant velocities $c_1$ and $c_2$ and their evolutions are controlled by two independent Poisson processes of rates 
$\lambda_1>0$ and $\lambda_2>0$, respectively, as described in Section 3 above. According to (\ref{gen4}), 
$e^{i\xi x_1^0} H_1(\xi,t)$ and $e^{i\xi x_2^0} H_2(\xi,t)$ are the characteristic functions of $X_1^{x_1^0}(t)$ and $X_2^{x_2^0}(t)$, respectively. 

Consider the sum $\tilde S(t) = X_1^{x_1^0}(t)+X_2^{x_2^0}(t)$ of these telegraph processes. The support of the distribution of $\tilde S(t)$ is 
the close interval $[(x_1^0+x_2^0)-(c_1+c_2)t, \; (x_1^0+x_2^0)+(c_1+c_2)t]$. 

If $c_1\neq c_2$ then the singular part of the distribution 
is concentrated at two terminal points $(x_1^0+x_2^0)\pm(c_1+c_2)t$ of this interval and  
$$\text{Pr} \{ \tilde S(t) = (x_1^0+x_2^0)\pm(c_1+c_2)t \} = \frac{1}{4} e^{-(\lambda_1+\lambda_2)t}, \qquad t> 0.$$ 
The density (in the sense of generalized functions) of the singular part of the distribution of $\tilde S(t)$ has the form 
$$\tilde\varphi_s(x,t) = \frac{e^{-(\lambda_1+\lambda_2)t}}{4} \left[ \delta((x_1^0+x_2^0)+(c_1+c_2)t) + \delta((x_1^0+x_2^0)-(c_1+c_2)t) \right],$$
where $\delta(x)$ is the Dirac delta-function. 
The absolutely continuous part of the distribution of $\tilde S(t)$ is concentrated in the open interval 
$((x_1^0+x_2^0)-(c_1+c_2)t, \; (x_1^0+x_2^0)+(c_1+c_2)t)$ and 
$$\text{Pr} \{ \tilde S(t) \in ((x_1^0+x_2^0)-(c_1+c_2)t, \; (x_1^0+x_2^0)+(c_1+c_2)t) \} = 1 - \frac{1}{2} e^{-(\lambda_1+\lambda_2)t} , \qquad t> 0.$$

If $c_1=c_2=c$ then the close interval $[(x_1^0+x_2^0)-2ct, \; (x_1^0+x_2^0)+2ct]$ is the support of the distribution of $\tilde S(t)$. 
The singular part of the distribution is concentrated at three points $x_1^0+x_2^0, \; (x_1^0+x_2^0)\pm 2ct$ of this interval and 
$$\text{Pr} \{ \tilde S(t) = (x_1^0+x_2^0)\pm 2ct \} = \frac{1}{4} e^{-(\lambda_1+\lambda_2)t}, \qquad t> 0.$$
$$\text{Pr} \{ \tilde S(t) = x_1^0+x_2^0 \} = \frac{1}{2} e^{-(\lambda_1+\lambda_2)t}, \qquad t> 0.$$ 
The density (in the sense of generalized functions) of the singular part of the distribution of $\tilde S(t)$ in this case has the form 
$$\tilde\varphi_s(x,t) = \frac{e^{-(\lambda_1+\lambda_2)t}}{2} \delta(x_1^0+x_2^0) + \frac{e^{-(\lambda_1+\lambda_2)t}}{4} \left[ \delta((x_1^0+x_2^0)+2ct) + \delta((x_1^0+x_2^0)-2ct) \right].$$
The absolutely continuous part of the distribution of $\tilde S(t)$ is concentrated in the area  
$\tilde M_t=((x_1^0+x_2^0)-2ct, \; x_1^0+x_2^0)\cup (x_1^0+x_2^0, \;   (x_1^0+x_2^0)+2ct)$ and 
$$\text{Pr} \{ \tilde S(t) \in \tilde M_t \} = 1 - e^{-(\lambda_1+\lambda_2)t} , \qquad t> 0.$$

The characteristic function of process $\tilde S(t)$ is given by 
$$\tilde\Psi(\xi,t) = e^{i\xi (x_1^0+x_2^0)} H_1(\xi,t) H_2(\xi,t), \qquad \xi\in\Bbb R, \quad t\ge 0.$$ 
If the start points $x_1^0, \;x_2^0$ are symmetric with respect to the origin $x=0$, then 
$x_1^0+x_2^0=0$ and in this case $\tilde \Psi(\xi,t)$ is a real-valued function, otherwise it is a complex function. Clearly, 
$\tilde\Psi(\xi,t)$ has a much more complicated form (in comparison with characteristic function $\Psi(\xi,t)$ given by (\ref{sum3})) that 
substantially depends on the numbers $\lambda_1/c_1$ and $\lambda_2/c_2$. 

To obtain the distribution of process $\tilde S(t)$ one needs to evaluate the inverse Fourier transform of the characteristic function 
$\tilde\Psi(\xi,t)$, however this is a very difficult problem that can, apparently, be done numerically only.

\section{Appendix}

\numberwithin{equation}{section}

In this appendix we prove four auxiliary lemmas that have been used in our analysis.

\bigskip

{\bf Lemma A1.} {\it For arbitrary positive $a>0, \; b>0$ the following formula holds:}
\begin{equation}\label{lem1}
\int_{-a}^a I_0\left( b\sqrt{a^2-x^2} \right) dx = \frac{2}{b} \sinh(ab) , \qquad a>0, \quad b>0. 
\end{equation}

\vskip 0.2cm

{\it Proof.} Using series representation (\ref{pprop5}) of the modified Bessel function $I_0(z)$, we get: 
$$\aligned 
\int_{-a}^a I_0\left( b\sqrt{a^2-x^2} \right) dx & = 2 \int_0^a I_0\left( b\sqrt{a^2-x^2} \right) dx \\
& = 2 \sum_{k=0}^{\infty} \frac{1}{(k!)^2} \left( \frac{b}{2} \right)^{2k} \int_0^a (a^2-x^2)^k \; dx \\
& = 2 \sum_{k=0}^{\infty} \frac{1}{(k!)^2} \left( \frac{b}{2} \right)^{2k} a^{2k+1} \int_0^1 (1-z^2)^k \; dz \\
& = 2 \sum_{k=0}^{\infty} \frac{1}{(k!)^2} \left( \frac{b}{2} \right)^{2k} a^{2k+1} \; \frac{\sqrt{\pi} \;\; k!}{(2k+1) \; \Gamma\left( k+\frac{1}{2} \right)} \\
& = 2 \sum_{k=0}^{\infty} \frac{1}{k!} \left( \frac{b}{2} \right)^{2k} a^{2k+1} \; \frac{\sqrt{\pi} \;\; 2^k}{(2k+1) \; \sqrt{\pi} \; (2k-1)!!} \\
& = 2 \sum_{k=0}^{\infty} \frac{a^{2k+1} \; b^{2k}}{(2k)!! \;\; (2k+1)!!} \\
& = \frac{2}{b} \sum_{k=0}^{\infty} \frac{(ab)^{2k+1}}{(2k+1)!} \\
& = \frac{2}{b} \sinh(ab) ,
\endaligned$$
where we have used the well-known formulas
$$\Gamma(z+1)=z\Gamma(z), \quad \Gamma\left( k+\frac{1}{2} \right) = \frac{\sqrt{\pi}}{2^k} (2k-1)!! , \quad (2k)!!=2^k \; k!,  \qquad k\ge 0, \;\; (-1)!!=1 .$$
The lemma is proved. $\square$ 

\bigskip

In particular, for $a=2ct, \; b=\lambda/c$ formula (\ref{lem1}) yields: 
\begin{equation}\label{lemm1}
\int_{-2ct}^{2ct} I_0\left( \frac{\lambda}{c} \sqrt{4c^2t^2-x^2} \right) dx = \frac{2c}{\lambda} \sinh(2\lambda t) . 
\end{equation}

The next two lemmas deal with the Fourier transformation
$$\mathcal F_x[ f(x) ](\xi) \equiv \hat f(\xi) = \int_{-\infty}^{\infty} e^{i\xi x} \; f(x) \; dx , \qquad \xi \in\Bbb R ,$$
and the inverse Fourier transformation
$$\mathcal F_{\xi}^{-1} [ \hat f(\xi) ](x) = \int_{-\infty}^{\infty} e^{-i\xi x} \; \hat f(\xi) \; d\xi , \qquad x\in\Bbb R ,$$
of the modified Bessel function $I_0(z)$.

\bigskip

{\bf Lemma A2.} {\it For arbitrary positive $a>0, \; b>0$ the following formula holds:}
\begin{equation}\label{lem2}
\mathcal F_x \left[ I_0(b\sqrt{a^2-x^2}) \; \Theta(a-|x|) \right](\xi) =  2 \left[ \frac{\sinh(a\sqrt{b^2-\xi^2})}{\sqrt{b^2-\xi^2}} \; \bold 1_{\{ |\xi| \le b \}}  + \frac{\sin(a\sqrt{\xi^2-b^2})}{\sqrt{\xi^2-b^2}} \; \bold 1_{\{ |\xi| > b  \}} \right] ,
\end{equation}
{\it where $\Theta(x)$ is the Heaviside step function defined by} (\ref{pprop4}) {\it and $\bold 1_{\{ z \}}$ is the indicator function}.

\vskip 0.2cm

{\it Proof.} We have
$$\aligned
\mathcal F_x \left[ I_0(b\sqrt{a^2-x^2}) \; \Theta(a-|x|) \right](\xi) & = \int\limits_{-a}^a e^{i\xi x} \; I_0(b\sqrt{a^2-x^2}) \; dx \\
& = 2 \int_0^a \cos{(\xi x)} \; I_0(b\sqrt{a^2-x^2}) \; dx \\
& \quad \text{(substitution $z=\sqrt{a^2-x^2}$)} \\
& = 2 \int_0^a \frac{z \; \cos{(\xi \sqrt{a^2-z^2})}}{\sqrt{a^2-z^2}} \; I_0(bz) \; dz \\
& \quad \text{(see \cite[item 2.15.10, formula 8]{pbm1})} \\
& = 2 \; \frac{\sin(a\sqrt{\xi^2-b^2})}{\sqrt{\xi^2-b^2}} \\
& = 2 \left[ \frac{\sinh(a\sqrt{b^2-\xi^2})}{\sqrt{b^2-\xi^2}} \; \bold 1_{\{ |\xi| \le b \}}  + \frac{\sin(a\sqrt{\xi^2-b^2})}{\sqrt{\xi^2-b^2}} \;
\bold 1_{\{ |\xi| > b  \}} \right] .
\endaligned$$
The lemma is proved. $\square$

\bigskip

In particular, by setting $a=2ct$ and $b=\lambda/c$ in (\ref{lem2}) (for arbitrary $c>0, \; \lambda>0, \; t>0$), we derive the following equality:
\begin{equation}\label{lem3}
\aligned
\mathcal F_x & \left[ I_0 \left( \frac{\lambda}{c} \sqrt{4c^2t^2-x^2} \right) \Theta(2ct-|x|) \right](\xi) \\
& \qquad = 2c \left[ \frac{\sinh(2t\sqrt{\lambda^2-c^2\xi^2})}{\sqrt{\lambda^2-c^2\xi^2}} \; \bold 1_{\{ |\xi| \le\frac{\lambda}{c} \}}  + \frac{\sin(2t\sqrt{c^2\xi^2-\lambda^2})}{\sqrt{c^2\xi^2-\lambda^2}} \; \bold 1_{\{ |\xi| > \frac{\lambda}{c} \}} \right] .
\endaligned
\end{equation}
Differentiating (\ref{lem3}) in $t$ we obtain
\begin{equation}\label{lem4}
\aligned
\mathcal F_x & \left[ \frac{\partial}{\partial t} I_0 \left( \frac{\lambda}{c} \sqrt{4c^2t^2-x^2} \right) \Theta(2ct-|x|) \right](\xi) + 4c \cos(2ct\xi) \\
& \qquad = 4c \left[ \cosh(2t\sqrt{\lambda^2-c^2\xi^2}) \; \bold 1_{\{ |\xi| \le\frac{\lambda}{c} \}}  + \cos(2t\sqrt{c^2\xi^2-\lambda^2}) \; \bold 1_{\{ |\xi| > \frac{\lambda}{c} \}} \right] .
\endaligned
\end{equation}
Applying inverse Fourier transformation to (\ref{lem3}) and (\ref{lem4}) we obtain
\begin{equation}\label{lem5}
\aligned
\mathcal F_{\xi}^{-1} \biggl[ \frac{\sinh(2t\sqrt{\lambda^2-c^2\xi^2})}{\sqrt{\lambda^2-c^2\xi^2}} \; \bold 1_{\{ |\xi| \le\frac{\lambda}{c} \}} & + \frac{\sin(2t\sqrt{c^2\xi^2-\lambda^2})}{\sqrt{c^2\xi^2-\lambda^2}} \; \bold 1_{\{ |\xi| > \frac{\lambda}{c} \}} \biggr](x) \\
& = \frac{1}{2c} \; I_0 \left( \frac{\lambda}{c} \sqrt{4c^2t^2-x^2} \right) \Theta(2ct-|x|) ,
\endaligned
\end{equation}
\begin{equation}\label{lem6}
\aligned
\mathcal F_{\xi}^{-1} \biggl[ & \cosh(2t\sqrt{\lambda^2-c^2\xi^2}) \; \bold 1_{\{ |\xi| \le\frac{\lambda}{c} \}} + \cos(2t\sqrt{c^2\xi^2-\lambda^2}) \; \bold 1_{\{ |\xi| > \frac{\lambda}{c} \}} \biggr] (x) \\
& = \frac{1}{4c} \; \frac{\partial}{\partial t} I_0 \left( \frac{\lambda}{c} \sqrt{4c^2t^2-x^2} \right) \Theta(2ct-|x|) + \frac{1}{2} \bigl[ \delta(2ct-x) + \delta(2ct+x) \bigr] ,
\endaligned
\end{equation}
where $\delta(x)$ is the Dirac delta function.

\bigskip

{\bf Lemma A3.} {\it For arbitrary positive $p>0, \; q>0,$ the following formula holds:}
\begin{equation}\label{lem7}
\aligned
\mathcal F_{\xi}^{-1} & \biggl[ \frac{\sinh^2(q\sqrt{p^2-\xi^2})}{p^2-\xi^2} \; \bold 1_{\{ |\xi| \le p \}} + \frac{\sin^2(q\sqrt{\xi^2-p^2})}{\xi^2-p^2} \; \bold 1_{\{ |\xi| > p \}} \biggr](x) \\
& \qquad\qquad = \frac{1}{4} \left\{ \int_{|x|}^{2q} I_0 \left( p \sqrt{\tau^2-x^2} \right) d\tau \right\}  \Theta(2q-|x|).
\endaligned
\end{equation}

\vskip 0.2cm

{\it Proof.} Applying Fourier transformation to the right-hand side of (\ref{lem7}) and using formula (\ref{lem2}), we have: 
$$\aligned
\frac{1}{4} & \mathcal F_x \left[ \left\{ \int_{|x|}^{2q} I_0 \left( p \sqrt{\tau^2-x^2} \right) d\tau \right\} \Theta(2q-|x|) \right](\xi) \\
& = \frac{1}{4} \mathcal F_x \left[ \left\{ \int_0^{2q} I_0 \left( p \sqrt{\tau^2-x^2} \right) \Theta(\tau-|x|) \; d\tau \right\} \Theta(2q-|x|) \right](\xi) \\ 
& = \frac{1}{4} \int_{-2q}^{2q} e^{i\xi x} \left\{ \int_0^{2q} I_0 \left( p \sqrt{\tau^2-x^2} \right) \Theta(\tau-|x|) \; d\tau \right\} dx \\
& = \frac{1}{4} \int_0^{2q} \left\{ \int_{-2q}^{2q} e^{i\xi x} \; I_0 \left( p \sqrt{\tau^2-x^2} \right) \Theta(\tau-|x|) \; dx \right\} d\tau \\
& = \frac{1}{4} \int_0^{2q} \left\{ \int_{-\tau}^{\tau} e^{i\xi x} \; I_0 \left( p \sqrt{\tau^2-x^2} \right) \; dx \right\} d\tau \\
& = \frac{1}{4} \int_0^{2q} \left\{ \mathcal F_x \left[ I_0(p\sqrt{\tau^2-x^2}) \; \Theta(\tau-|x|) \right](\xi) \right\} d\tau \\
& = \frac{1}{2} \int_0^{2q} \left[ \frac{\sinh(\tau\sqrt{p^2-\xi^2})}{\sqrt{p^2-\xi^2}} \bold 1_{\{ |\xi| \le p \}} +  \frac{\sin(\tau\sqrt{\xi^2-p^2})}{\sqrt{\xi^2-p^2}} \bold 1_{\{ |\xi| > p \}} \right] d\tau \\
& = \frac{1}{2\sqrt{p^2-\xi^2}} \left\{ \int_0^{2q} \sinh(\tau\sqrt{p^2-\xi^2}) \; d\tau \right\} \bold 1_{\{ |\xi| \le p \}} \\
& \hskip 2cm + \frac{1}{2\sqrt{\xi^2-p^2}} \left\{ \int_0^{2q} \sin(\tau\sqrt{\xi^2-p^2}) \; d\tau \right\} \bold 1_{\{ |\xi| > p \}} \\
& = \frac{1}{p^2-\xi^2} \; \frac{\cosh(2q\sqrt{p^2-\xi^2}) - 1}{2} \; \bold 1_{\{ |\xi| \le p \}} +  \frac{1}{\xi^2-p^2} \; \frac{1 - \cos(2q\sqrt{\xi^2-p^2})}{2} \; \bold 1_{\{ |\xi| > p \}} \\
& = \frac{\sinh^2(q\sqrt{p^2-\xi^2})}{p^2-\xi^2} \; \bold 1_{\{ |\xi| \le p \}} + \frac{\sin^2(q\sqrt{\xi^2-p^2})}{\xi^2-p^2} \; \bold 1_{\{ |\xi| > p \}} .
\endaligned$$
The lemma is proved. $\square$

\bigskip 

In particular, setting $q=ct, \; p=\lambda/c$ in (\ref{lem7}) we arrive to the formula 
\begin{equation}\label{lem8}
\aligned
\mathcal F_{\xi}^{-1} & \biggl[ \frac{\sinh^2\left( t\sqrt{\lambda^2-c^2\xi^2}\right) }{\lambda^2-c^2\xi^2} \; \bold 1_{\{ |\xi| \le \frac{\lambda}{c} \}} + \frac{\sin^2\left( t\sqrt{c^2\xi^2-\lambda^2}\right) }{c^2\xi^2-\lambda^2} \; \bold 1_{\{ |\xi| > \frac{\lambda}{c} \}} \biggr](x) \\
& \qquad\qquad = \frac{1}{4c^2} \left\{ \int_{|x|}^{2ct} I_0 \left( \frac{\lambda}{c} \sqrt{\tau^2-x^2} \right) d\tau \right\} \Theta(2ct-|x|).
\endaligned
\end{equation}
Note that differentiating (\ref{lem8}) in $t$, we obtain again formula (\ref{lem5}).

\bigskip

{\bf Lemma A4.} {\it For arbitrary integers $n\ge 0, \; k\ge 0$ such that $n\ge 2k$ and for arbitrary real $x\in\Bbb R$ the following formula holds:}
\begin{equation}\label{hyp1}
\int z^n F\left( -k, \frac{1}{2}; \; \frac{3}{2};  \; \frac{x^2}{z^2} \right) dz = \frac{z^{n+1}}{n+1} \; _3F_2 \left( -k, -\frac{n}{2}-\frac{1}{2},  \frac{1}{2}; \; -\frac{n}{2}+\frac{1}{2}, \frac{3}{2}; \; \frac{x^2}{z^2} \right) + \psi, 
\end{equation}
{\it where the hypergeometric function on the right-hand side of} (\ref{hyp1}) {\it is defined by} (\ref{dist3}) {\it and $\psi$ is an arbitrary function not depending on $z$.} 

\vskip 0.2cm

{\it Proof.} Differentiating in $z$ the function on the right-hand side of (\ref{hyp1}) and using the second formula of (\ref{poch}) we obtain: 
$$\aligned 
\frac{1}{n+1} & \; \frac{d}{dz} \left[ z^{n+1} \; _3F_2 \left( -k, -\frac{n}{2}-\frac{1}{2},  \frac{1}{2}; \; -\frac{n}{2}+\frac{1}{2}, \frac{3}{2}; \; 
\frac{x^2}{z^2} \right) \right] \\ 
& = \frac{1}{n+1} \; \frac{d}{dz} \sum_{s=0}^k \frac{(-k)_s \; \left( \frac{1}{2} \right)_s}{\left( \frac{3}{2} \right)_s} \; \frac{\left( -\frac{n}{2}-\frac{1}{2} \right)_s}{\left( \left(-\frac{n}{2}-\frac{1}{2}\right) + 1 \right)_s} \; \frac{x^{2s}}{s!} \; z^{n-2s+1} \\ 
& = \frac{1}{n+1} \; \frac{d}{dz} \sum_{s=0}^k \frac{(-k)_s \; \left( \frac{1}{2} \right)_s}{\left( \frac{3}{2} \right)_s} \; \frac{n+1}{n-2s+1} \; \frac{x^{2s}}{s!} \; z^{n-2s+1} \\ 
& = z^n \sum_{s=0}^k \frac{(-k)_s \; \left( \frac{1}{2} \right)_s}{\left( \frac{3}{2} \right)_s} \; \frac{1}{s!} \; \left(\frac{x^2}{z^2}\right)^s \\ 
& = z^n F\left( -k, \frac{1}{2}; \; \frac{3}{2};  \; \frac{x^2}{z^2} \right) , 
\endaligned$$
coinciding with integrand on the left-hand side of (\ref{hyp1}). The lemma is proved. $\square$


\begin{thebibliography}{29}

\bibitem{bart1}
Bartlett M. Some problems associated with random velocity. {\it
Publ. Inst. Stat. Univ. Paris,} 1957, {\bf 6}, 261-270.

\bibitem{bart2}
Bartlett M. A note on random walks at constant speed. {\it Adv.
Appl. Probab.,} 1978, {\bf 10}, 704-707.

\bibitem{br}
Bogachev L., Ratanov N. Occupation time distributions for the
telegraph process. {\it Stoch. Process. Appl.,} 2011, {\bf 121},
1816-1844.

\bibitem{cane1}
Cane V. Random walks and physical processes. {\it Bull. Intern.
Statist. Inst.,} 1967, {\bf 42}, 622-640.

\bibitem{cane2}
Cane V. Diffusion models with relativity effects. // In: {\it
Perspectives in Probability and Statistics,} Sheffield, Applied
Probability Trust, 1975, 263-273.

\bibitem{cres1}
Di Crescenzo A. On random motion with velocities alternating at
Erlang-distributed random times. {\it Adv. Appl. Probab.,} 2001,
{\bf 33}, 690-701.

\bibitem{cres2}
Di Crescenzo A., Martinucci B. A damped telegraph random process
with logistic stationary distributions. {\it J. Appl. Probab.,}
2010, {\bf 47}, 84-96.

\bibitem{foong1}
Foong S.K. First-passage time, maximum displacement and Kac's
solution of the telegrapher's equation. {\it Phys. Rev. A,} 1992,
{\bf 46}, 707-710.

\bibitem{foong2}
Foong S.K., Kanno S. Properties of the telegrapher's random
process with or without a trap. {\it Stoch. Process. Appl.,} 2002,
{\bf 53}, 147-173.

\bibitem{gold}
Goldstein S. On diffusion by discontinuous movements and on the
telegraph equation. {\it Quart. J. Mech. Appl. Math.}, 1951, {\bf
4}, 129-156.

\bibitem{kab}
Kabanov Yu.M. Probabilistic representation of a solution of the
telegraph equation. {\it Theory Probab. Appl.}, 1992, {\bf 37},
379-380.

\bibitem{kac}
Kac M. A stochastic model related to the telegrapher's equation.
{\it Rocky Mount. J. Math.}, 1974, {\bf 4}, 497-509.

\bibitem{kap}
Kaplan S. Differential equations in which the Poisson process
plays a role. {\it Bull. Amer. Math. Soc.,} 1964, {\bf 70},
264-267.

\bibitem{kis}
Kisynski J. On M.Kac's probabilistic formula for the solution of
the telegraphist's equation. {\it Ann. Polon. Math.}, 1974, {\bf
29}, 259-272.

\bibitem{kol1}
Kolesnik A.D. The equations of Markovian random evolution on the
line. {\it J. Appl. Probab.,} 1998, {\bf 35}, 27-35.

\bibitem{kol2}
Kolesnik A.D. Moment analysis of the telegraph random process. 
{\it Bull. Acad. Sci. Moldova, Ser. Math.}, 2012, {\bf 1(68)}, 90-107.

\bibitem{kol3}
Kolesnik A.D. Probability distribution function for the Euclidean distance
between two telegraph processes. {\it Adv. Appl. Probab.,} 2014, {\bf 46}. 
(To appear, electronic preprint arXiv:1305.6522) 

\bibitem{kol4}
Kolesnik A.D., Ratanov N. {\it Telegraph Processes and Option Pricing}. 
Springer, 2013, Heidelberg. 

\bibitem{mas1}
Masoliver J., Weiss G.H. First-passage times for a generalized
telegrapher's equation. {\it Physica A}, 1992, {\bf 183}, 537-548.

\bibitem{mas2}
Masoliver J., Weiss G.H. On the maximum displacement of a
one-dimensional diffusion process described by the telegrapher's
equation. {\it Physica A}, 1993, {\bf 195}, 93-100.

\bibitem{pin}
Pinsky M.A. {\it Lectures on Random Evolution}. World Sci., 1991,
River Edge, NJ.

\bibitem{pbm1}
Prudnikov A.P., Brychkov Yu.A., Marichev O.I. {\it Integrals and
Series. Special Functions.} Nauka, 1983, Moscow. (In Russian)

\bibitem{pbm2}
Prudnikov A.P., Brychkov Yu.A., Marichev O.I. {\it Integrals and
Series. Additional Chapters.} Nauka, 1986, Moscow. (In Russian)

\bibitem{rat1}
Ratanov, N. Random walks in an inhomogeneous one-dimensional
medium with reflecting and absorbing barriers. 
{\it Theoret. Math. Phys.,} 1997, {\bf 112}, 857-865.

\bibitem{rat2}
Ratanov, N. Telegraph evolutions in inhomogeneous media. {\it Markov
Process. Related Fields}, 1999, {\bf 5}, 53-68.

\bibitem{sta}
Stadje W., Zacks S. Telegraph processes with random velocities.
{\it J. Appl. Probab.,} 2004, {\bf 41}, 665-678.

\bibitem{turb}
Turbin A.F., Samoilenko I.V. A probabilistic method for solving
the telegraph equation with real-analytic initial conditions.
{\it Ukrain. Math. J.}, 2000, {\bf 52}, 1292-1299.

\end{thebibliography}
\end{document}